\documentclass[12pt]{amsart}

\usepackage{amsfonts}
\usepackage{pstricks, pst-node, pst-text, pst-3d}
\usepackage{epsfig}
\usepackage{graphicx}
\usepackage{amsmath}
\usepackage{amssymb}

 \advance\hoffset by -0.0 truecm

\newtheorem{theorem}{Theorem}

\newtheorem{proposition}{Proposition}
\theoremstyle{definition}

\theoremstyle{remark}
\newtheorem{remark}{Remark}
\numberwithin{equation}{section}

\newcommand{\E}{\mathcal{E}}

\newcommand{\re}{{\rm Re\, }}

\newcommand{\bu}{{\bf u}}

\DeclareMathOperator{\spn}{span}

\DeclareMathOperator{\cotan}{cotan}

\begin{document}

\title[Modified action and differential operators]
{Modified action and differential operators on the 3-D sub-Riemannian sphere}

\author{Der-Chen Chang,\ Irina Markina, and Alexander Vasil'ev}

\address{Department of Mathematics, Georgetown University, Washington
D.C. 20057, USA}

\email{chang@georgetown.edu}

\thanks{}

\address{Department of Mathematics,
University of Bergen, Johannes Brunsgate 12, Bergen 5008, Norway}

\email{irina.markina@uib.no} \email{alexander.vasiliev@uib.no}

\thanks{The first author is partially supported by a research
grant from the United State Army Research Office and a Hong Kong RGC
competitive earmarked research grant $\#$600607. The second and the
third authors have been  supported by the grant of the Norwegian Research Council \#177355/V30, and by the European Science Foundation Research Networking Programme HCAA}




\date{July 26, 2008}


\begin{abstract}
Our main aim is to present a geometrically meaningful formula for the fundamental solutions to a second order sub-elliptic differential equation and to the  heat equation associated with a sub-elliptic operator in the sub-Riemannian geometry on the unit sphere $\mathbb S^3$.
Our method is based on the Hamiltonian approach, where the corresponding Hamitonian system is solved with mixed boundary conditions. A closed form of the  modified action is given. It is a sub-Riemannian invariant and plays the role of a distance on $\mathbb S^3$. 
\end{abstract}

\maketitle

\section{Introduction}

The unit $3$-sphere centered on the origin is a subset of $\mathbb
R^4$ defined as
$$\mathbb S^3=\{(x_1,x_2,x_3,x_4)\in\mathbb R^4:\ x_1^2+x_2^2+x_3^2+x_4^2=1)\}.$$
Regarding $\mathbb R^4$ as the space of
quaternions $\mathbb H$, the above unit 3-sphere admits the form
$$\mathbb S^3=\{q\in\mathbb H:\ |q|^2=1)\}.$$ 
This description represents the sphere $\mathbb S^3$ as a set of
unit quaternions with the inherited group structure, and it can be considered as the spin group $Sp(1)$,
where the group operation is just the multiplication of quaternions. Let us identify $\mathbb R^3$ with
pure imaginary quaternions. The conjugation $qh\bar q$ of a pure
imaginary quaternion $h$ with a unit quaternion $q$ defines
rotation in $\mathbb R^3$, and since $|qh\bar q|=|h|$, the map
$h\mapsto qh\bar q$ defines a two-to-one homomorphism $Sp(1)\to
SO(3)$. The Hopf map $\pi:\mathbb S^3\to\mathbb S^2$ can be
defined by $$\mathbb S^3\ni q\mapsto qi\bar q=\pi(q)\in\mathbb
S^2.$$ In its turn, the Hopf map defines a principle circle bundle also known
as the Hopf bundle.

The sub-Riemannian structure of $\mathbb S^3$ comes naturally from the non-commutative
group structure of the sphere in the sense that two vector fields
span the smoothly varying distribution of the tangent bundle, and
their commutator generates the missing direction. The missing direction is also can be obtained as an integral line of the Hopf vector field corresponding to the Hopf fibration. The sub-Riemannian geometry on $\mathbb S^3$ was studied in~\cite{CalinChangMarkina1, ChangMarkinaVasiliev, HurtadoRosales}, see also \cite{Boscain2}. Explicit formulas for  geodesics were given in \cite{ChangMarkinaVasiliev}. Let us mention that the word `geodesic' in our terminology stands for the projection of
the solutions to a Hamiltonian system onto the underlying manifold, that is a good generalization of the notion of  geodesic from Riemannian to  sub-Riemannian manifolds, see for instance~\cite{Montgomery, Strichartz}.
The Lagrangian approach was applied in~\cite{CalinChangMarkina1} and~\cite{HurtadoRosales} in order to characterize and to find the shortest geodesics. Another approach based on the control theory was employed in~\cite{Boscain2}.  

In this paper, our main aim is to deduce a geometrically meaningful formula for the Green function for a second order sub-elliptic differential operator and the heat kernel associated with this operator in the sub-Riemannian geometry on the unit sphere $\mathbb S^3$. There exists a vast amount of  literature studying sub-elliptic operators based on different methods. Here we give only  few  possible references~\cite{FollandStein, Geller, Hormander, KohnRossi, RothschildStein, Taylor}. The exact form of the heat kernel on sub-Riemannian manifolds was obtained only in some
simple particular cases. Namely, in the case of the Heisenberg group, the representation theory was used in \cite{Hulanicki}, the probability approach
was developed in \cite{Gaveau}, the Laguerre calculus was applied in \cite{BealsGreiner}, and the Hamilton-Jacobi approach one finds in \cite{BealsGaveauGreiner3}. Recently, the heat kernel was investigated on other low dimensional sub-Riemannian manifolds. The representation
theory was exploiting in  \cite{Agrachev} and \cite{Bauer}, spectral analysis and small time asymptotics were used in \cite{Baudoin}.
Our method is based on the Hamiltonian approach which benefit is the connection between the heat kernel and the geometry of the sphere
as a sub-Riemannian manifold.
Analogously to Hadamard's method for strictly hyperbolic operators, our method essentially uses three important ingredients:
\begin{itemize}
\item Solution of the Hamiltonian system with non-standard boundary conditions and construction of a modified action on solutions to this systems. This modified action plays the role of a sub-Riemannian distance;

\item Solution of the corresponding transport equations and deduction of the volume elements.
\item Integration of the  modified action over the characteristic variety with respect the measure defined by the volume element.
\end{itemize}

This method was realized for the two step nilpotent groups, for instance, in series of papers~\cite{BealsGaveauGreiner1, BealsGaveauGreiner2, BealsGaveauGreiner3}, where the geometric meaning of the fundamental solutions was revealed. For other geometries see, for example~\cite{BealsGaveauGreinerKannai, BealsGaveauGreiner4, CalinChangMarkina2}. The case of 3-sphere reveals new features, and possessing the Cartan decomposition of the acting group,  it is not a direct analog of previous considerations.

The structure of the paper is as follows. The classical setup for the heat kernel in the Riemannian case is presented in Section 2.  In Section 3, we define the horizontal distribution and the sub-Riemannian metric. The Hamiltonian system is derived in the fourth section. In the fifth section we treat the problem of finding geodesics as an optimal control problem. Symmetries of the Hamiltonain system are discussed. In  Section 6, we solve the Hamiltonian system  to find geodesics and to solve the boundary value problem. The number of geodesics connecting two fixed points on $\mathbb S^3$ is studied. Both cartesian and hyperspherical coordinates are used. At the end of this section we define the modified action and investigate its properties.  Special directions in the cotangent bundle given by the Hamiltonian system are revealed clearly in the hyperspherical coordinates contrasting with the cartesian ones. We use these directions to construct the modified action solving  the Hamiltonian system with non-standard mixed boundary conditions. The modified action satisfies a generalized Hamilton-Jacobi equation (Section 7). It is a sub-Riemannian invariant on $\mathbb S^3$ and it is used for the construction of a distance function (Section 8). The distance function is involved into the fundamental solutions to the sub-Laplacian equation and to the heat equation associated with the sub-Laplacian.  The concluding Section 9 is concerned with the volume element. The sub-Laplacian and the heat operator associated with this sub-Laplacian are not elliptic, they  degenerate along a singular manifold of dimension one in the cotangent space. The fundamental solutions to these equations can be obtained by integrating the distance function over this one-dimensional singular set which is the characteristic variety of the corresponding Hamiltonian  with respect of a special measure with the density called the {\it volume element}. Unlike the case of  nilpotent groups the volume element depends on  phase variables that does not permit to find its explicit form. Instead we present differential equations, called the {\it transport equations} which solutions give the necessary volume elements.

The paper was initiated when the authors visited the National Center for
Theoretical Sciences and National Tsing Hua University during
May 2008. They would like to express their profound gratitude
to Professor Jing Yu  for the invitation and
for the warm hospitality of the staff extended to them during their stay in
Taiwan.

\section{Heat kernel in $\mathbb R^n$}

Let us present some simple calculations in $\mathbb R^n$ for the heat operator  motivating  further generalizations to the case of sub-Riemannian geometry on $\mathbb S^3$. Let $\Delta=\frac{1}{2}\sum_{j=1}^{n}\left(\frac{\partial}{\partial x_j}\right)^2$ be the Laplace operator. Then the kernel $P_u(x,x_0)$ for the operator $\Delta-\frac{\partial}{\partial u}$ is given by $$P_u(x,x_0)=\frac{1}{(2\pi u)^{\frac{n}{2}}}e^{-\frac{|x-x_0|^2}{2u}}.$$
If we write $f=\frac{1}{2}|x-x_0|^2$, then it is easy to see that the function $\frac{f}{u}$ satisfies the Hamilton-Jacobi equation $$\frac{\partial}{\partial u}\Big(\frac{f}{u}\Big)+\frac{1}{2}\sum_{j=1}^{n}\Big(\frac{\partial}{\partial x_j}\Big(\frac{f}{u}\Big)\Big)^2=0,\quad\mbox{with \ } \frac{1}{2}\sum_{j=1}^{n}\Big(\frac{\partial}{\partial x_j}\Big(\frac{f}{u}\Big)\Big)^2=H\Big(\nabla\big(\frac{f}{u}\big)\Big),$$ and $H$ is the Hamiltonian function associated with the Laplace operator $\Delta$. In the standard theory, the function $S=\frac{f}{u}$ is the classical action related to the Hamiltonian~$H$.

In the case of a general second order elliptic operator defined by smooth linearly independent vector fields $X_j$, $j=1,\ldots,n$ in $\mathbb R^n$, the heat kernel $P_u(x,x_0)$ for the operator $$\Delta_X-\frac{\partial}{\partial u}, \quad\text{with}\quad \Delta_X=\frac{1}{2}\sum_{j=1}^{n}X_j^2,$$ admits the form $$P_u(x,x_0)=\frac{1}{(2\pi u)^{\frac{n}{2}}}e^{-\frac{|x-x_0|^2}{2u}}(v_0+v_1u+v_2u^2+\ldots),$$ where the function $\frac{|x-x_0|^2}{2u}$ still satisfies the Hamilton-Jacoby equation with respect to the vector fields $X_j$. Associated Hamiltonian is degenerating only at one point  of $\mathbb R^n\times \mathbb R^n$ and the constants $v_l$ are chosen so that the delta function supported at $x_0$ is clearly seen.

Let us consider the vector fields $X_1,\ldots,X_k$ satisfying the Chow-Rashevski{\u\i} (or bracket generating) condition \cite{Chow, Rashevsky} (see Section~3) on $n$-dimensional manifold $M$, $k<n$. In this case the operator $\Delta_X=\frac{1}{2}\sum_{j=1}^{k}X_j$ is sub-elliptic and degenerates over a set of positive measure. Previous studies (see, e.g., \cite{BealsGaveauGreiner1, BealsGaveauGreiner2, BealsGaveauGreiner3, BealsGaveauGreiner4, BealsGaveauGreinerKannai,  CalinChangMarkina2} ) show that it is reasonable to expect the heat kernel $P_u(x,x_0)$ for the operator associated with the sub-Laplacian $\Delta_X$ in the form $$P_u(x,x_0)=\frac{C}{ u^q}\int\limits_{chv(H)_{x_0}(\tau)}e^{-\frac{f(x,x_0,\tau)}{u}}v(x,u,\tau)\,d\tau.$$ Here $chv(H)_{x_0}$ is the characteristic variety of the Hamiltonian function at $x_0$ associated with the sub-Laplacian $\Delta_X$ defined by $$chv(H)=\big\{(x,\xi)\in T^*M:\ H(x,\xi)=0\big\}.$$ The characteristic variety represents the singular set of the sub-elliptic operator. The function $f(x,x_0,\tau)$ plays the role of  square of the distance between the points $x_0$ and $x$ on the manifold $M$ and satisfies the generalized Hamilton-Jacobi equation $$\tau
\frac{d f}{d\tau}+H(x,\nabla_xf)=f.$$ The function $f$ is {\it a modified action}  associated with the degenerating Hamiltonian. The term $v(x,\tau)$ is a suitable measure on the characteristic variety $chv(H)_{x_0}$ at $x_0$ making the integral convergent. It is called the {\it  volume element} and it can be found from a differential equation known as {\it the transport equation}.

The following sections will be devoted to the study of the Hamiltonian system, its solutions and the construction of the modified action function $f$ as a distance function in the heat kernel associated to the sub-elliptic operator on $\mathbb S^3$ .

\section{Horizontal distribution on $\mathbb S^3$}

Let us turn to the sub-Riemannian geometry on  $\mathbb S^3$.
In order to calculate left-invariant vector fields we use the definition of  $\mathbb S^3$ as a set of unit
quaternions equipped with the following non-commutative multiplication `$\circ$': if $x=(x_1,x_2,x_3,x_4)$ and
$y=(y_1,y_2,y_3,y_4)$, then
\begin{eqnarray}\label{grlaw2}
x\circ y=(x_1,x_2,x_3,x_4)\circ(y_1,y_2,y_3,y_4) & = & \Big((x_1y_1-x_2y_2-x_3y_3-x_4y_4),
\nonumber \\ & \ \ & (x_2y_1+x_1y_2-x_4y_3+x_3y_4),
\nonumber
\\& \ & (x_3y_1+x_4y_2+x_1y_3-x_2y_4),\\ & \ & (x_4y_1-x_3y_2+x_2y_3+x_1y_4)\Big).\nonumber
\end{eqnarray}
The rule~\eqref{grlaw2} gives us the left translation $L_x(y)$ of an element $y=(y_1,y_2,y_3,y_4)$ by an element $x=(x_1,x_2,x_3,x_4)$. The left-invariant basis vector fields are defined as $X(x)=(L_x(y))_*X(0)$, where $X(0)$ are the basis vectors at the unity of the group.
Calculating the action of $(L_x(y))_*$ in the basis of the unit vectors of $\mathbb R^4$ we obtain four left-invariant vector fields
\begin{eqnarray}\label{vf}
X_1(x) & = & x_1\partial_{x_1}+x_2\partial_{x_2}+x_3\partial_{x_3}+x_4\partial_{x_4},\nonumber \\
X_2(x) & = & -x_2\partial_{x_1}+x_1\partial_{x_2}+x_4\partial_{x_3}-x_3\partial_{x_4},\\
X_3(x) & = & -x_3\partial_{x_1}-x_4\partial_{x_2}+x_1\partial_{x_3}+x_2\partial_{x_4},\nonumber \\
X_4(x) & = & -x_4\partial_{x_1}+x_3\partial_{x_2}-x_2\partial_{x_3}+x_1\partial_{x_4}.\nonumber
\end{eqnarray} It is easy to see that the vector $X_1(x)$ is the unit normal to $\mathbb S^3$ at $x$ with respect to the usual inner product $\langle\cdot,\cdot\rangle$ in $\mathbb R^4$, hence, we denote  $X_1(x)$ by $N$. Moreover, the vector fields $X_2(x)$, $X_3(x)$, $X_4(x)$ form an orthonormal basis  of the tangent space $T_x\mathbb S^3$ with respect to $\langle\cdot,\cdot\rangle$ at any point $x\in\mathbb S^3$. Let us denote these vector fields by $$X_3=X,\quad X_4=Y,\quad X_2=Z.$$

The vector fields possess the following commutation relations
$$[X,Y]=XY-YX=2Z,\quad [Z,X]=2Y,\quad [Y,Z]=2X.$$ Let $\mathcal{D} = \spn
\{ X, Y\}$ be the distribution generated by the vector fields $X$
and $Y$. Since $[X, Y] =  2Z \notin \mathcal{D} $, it follows that
$\mathcal{D}$ is not involutive. The distribution $\mathcal{D}$
will be called  {\it horizontal}. Any curve on the sphere with the
velocity vector contained in the distribution $\mathcal{D}$ will
be called a {\it horizontal curve}. Since $T_x\mathbb
S^3=\spn\{X,Y,Z=1/2[X,Y]\}$, the distribution is bracket
generating at each point $x\in \mathbb S^3$, see~\cite{Chow, Rashevsky}. We define the metric on the distribution $\mathcal D$
as the restriction of the metric $\langle\cdot,\cdot\rangle$ to
$\mathcal D$, and the same notation $\langle\cdot,\cdot\rangle$ will be used. 
This metric coincides with the metric given by the Killing form on the Lie algebra 
$T_e\mathbb S^3$.
Finally, the manifold
$(\mathbb S^3,\mathcal D,\langle\cdot,\cdot\rangle)$ becomes a step two
sub-Riemannian manifold.

\begin{remark}\label{rem1}
Observe that the choice of the horizontal distribution is not unique.
The relations $[Z,X]=2Y$ and $[Y,Z]=2X$ imply possible choices $\mathcal{D} =
\spn \{ X, Z\}$ or $\mathcal{D} = \spn \{ Y, Z\}$. The geometries
defined by different horizontal distributions are cyclically
symmetric, so we restrict our attention to the distribution $\mathcal{D} = \spn
\{ X, Y\}$.
\end{remark}
\begin{remark}\label{rem2}
Let us define two rotations in the planes $(x_1,x_2)$ and  $(x_3,x_4)$ as 
$$x=(x_1,x_2,x_3,x_4)\quad\mapsto\quad R_{1_{\phi}}(x)=(x_1\cos\phi-x_2\sin\phi,x_1\sin\phi+x_2\cos\phi,x_3,x_4),$$
$$x=(x_1,x_2,x_3,x_4)\quad\mapsto\quad R_{2_{\phi}}(x)=(x_1,x_2,x_3\cos\phi-x_4\sin\phi,x_3\sin\phi+x_4\cos\phi).$$ It is easy to see that these transformations leave $\mathbb S^3$ invariant. The vector fields $X$ and $Y$ change under these rotations as follows.  Under the rotation $R_{1_{\phi}}$ we have $$X\quad\mapsto\quad\widetilde X=X\cos\phi+Y\sin\phi,\qquad Y\quad\mapsto\quad\widetilde Y=-X\sin\phi+Y\cos\phi,$$ and under the transformation $R_{2_{\phi}}$ we have
$$X\quad\mapsto\quad\widetilde X=X\cos\phi-Y\sin\phi,\qquad Y\quad\mapsto\quad\widetilde Y=X\sin\phi+Y\cos\phi.$$ Since $[\widetilde X,\widetilde Y]=[X,Y]$, we conclude that these transformations preserve the horizontal distribution. In both cases the sub-Laplacian is also invariant $\Delta_X=X^2+Y^2=\widetilde X^2+\widetilde Y^2$.
\end{remark}

We also can define the distribution as a kernel of the following one-form
\begin{equation}\label{forma}
\omega = -x_2 dx_1 + x_1 dx_2 + x_4 dx_3 - x_3 dx_4
\end{equation}
 on $\mathbb{R}^4$.
One can easily check that
$$
\omega(X) =0,\quad \omega(Y)=0,\quad \omega (Z)=1 \not=0,\quad \omega (N)=0.$$
Hence, the horizontal distribution $\mathcal{D}_x$ at $x\in \mathbb S^3$ can be
written as $\ker \omega_x \cap T_x \mathbb{S}^3$. The one-form $\omega$ has the following geometric meaning. It is the difference of two independent area forms $\alpha=-x_2 dx_1 + x_1 dx_2$ in $(x_1,x_2)$-plane and $\beta=-x_4 dx_3 + x_3 dx_4$ in $(x_3,x_4)$-plane.

Let $\gamma(s)=(x_1(s),x_2(s),x_3(s),x_4(s))$ be a curve on $\mathbb S^3$. Then the velocity vector, written in the left-invariant basis, is $$\dot\gamma(s)=a(s)X(\gamma(s))+b(s)Y(\gamma(s))+c(s)Z(\gamma(s)),$$ where
\begin{eqnarray}\label{coord}
a & = & \langle\dot\gamma, X\rangle=-x_3\dot x_1-x_4\dot x_2+x_1\dot x_3+x_2\dot x_4,\nonumber \\
b & = & \langle\dot\gamma, Y\rangle=-x_4\dot x_1+x_3\dot x_2-x_2\dot x_3+x_1\dot x_4, \\
c & = & \langle\dot\gamma, Z\rangle=-x_2\dot x_1 +x_1\dot x_2+x_4\dot x_3-x_3\dot x_4.\nonumber
\end{eqnarray}
The following proposition holds.

\begin{proposition}
\label{prop1} Let $\gamma (s) = (x_1(s), x_2(s), x_3(s),
x_4(s))$ be a curve on $\mathbb{S}^3$. The curve $\gamma$ is
horizontal, if and only if,
\begin{equation}\label{hc}
c=\langle \dot\gamma, Z \rangle =-x_2\dot x_1 +x_1\dot x_2+x_4\dot x_3-x_3\dot x_4=0.
\end{equation}
\end{proposition} If we take into account the geometric meaning of the one-form $\omega$, then we can reformulate Proposition~\ref{prop1} in the following way. Let us denote by $A$ the area swept by the projection of the horizontal curve $\gamma$ onto the $(x_1,x_2)$-plane and bounded by the straight line connecting its ends, and by $B$ we denote the analogous area swept by the projection of the horizontal curve onto the $(x_3,x_4)$-plane. 
\begin{proposition}
\label{prop11} Let $\gamma (s) = (x_1(s), x_2(s), x_3(s),
x_4(s))$ be a curve on $\mathbb{S}^3$ and let $A$, $B$ be as introduced above. Then, the curve $\gamma$ is
horizontal, if and only if, $A=B$.
\end{proposition}

\begin{figure}[ht]
\hspace{2.5cm}
\begin{pspicture}(2,1)(15,9)
\psline[linewidth=0.15mm]{-}(2,2)(6,4)
\psline[linewidth=0.15mm]{-}(2,2)(2,7)
\psline[linewidth=0.15mm]{-}(6,4)(6,9)
\psline[linewidth=0.15mm]{-}(2,7)(6,9)
\psline[linewidth=0.15mm]{-}(2,2)(9,2)
\psline[linewidth=0.15mm]{-}(6,4)(13,4)
\psline[linewidth=0.15mm]{-}(9,2)(13,4)
\pscurve[linewidth=0.8mm](4,3)(6,3.3)(7,3.5)(8,3.75)(10,4.7)(10.4,5.7)(10,6.6)(8,7)
 \pscircle[fillstyle=solid,
fillcolor=black](4,3){.1} 
 \pscircle[fillstyle=solid,
fillcolor=black](8,7){.1} 
\psline[linewidth=0.15mm, linestyle=dashed](8,7)(8,3.5)
\psline[linewidth=0.15mm, linestyle=dashed](8,7)(5,7)
\pscurve[fillstyle=solid,fillcolor=lightgray]
(4,3)(5,3)(6,3)(7,3)(8,3.03)(9,3.07)(10,3.13)(10.4,3.25)
(10,3.41)(9.5,3.46)(9,3.48)(8,3.5)
\pscurve[fillstyle=solid,fillcolor=lightgray]%
(4,3)(4.1,5)(4.4,6.2)(5,7)
\rput(3.6, 3.1){0}
\rput(3, 6.8){$(x_1,x_2)$}
\rput(8.3, 2.4){$(x_3,x_4)$}
\rput(10.7, 5.5){$\gamma$}
\rput(4.5, 5.7){$A$}
\rput(8, 3.3){$B$}
\end{pspicture}
\caption[]{Projections of $\gamma$ to the planes  $(x_1,x_2)$ and $(x_3,x_4)$ in Proposition~\ref{prop11} }\label{fig1}
\end{figure}
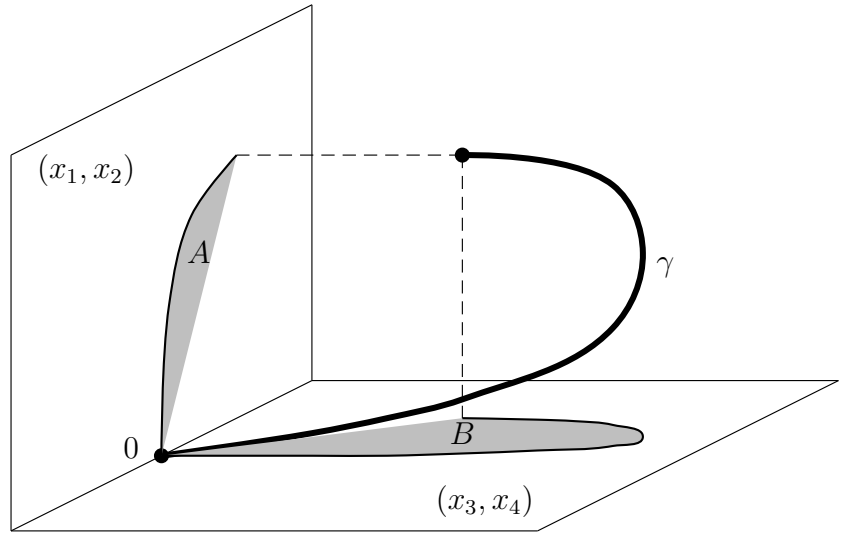

The manifold $\mathbb S^3$ is connected and it satisfies the
bracket generating condition. By the Chow-Rashevski{\u\i} theorem~\cite{Chow, Rashevsky},
there exist piecewise $C^1$ horizontal curves connecting two
arbitrary points of $\mathbb S^3$. In fact, smooth horizontal
curves connecting two arbitrary points of $\mathbb S^3$ were
constructed in~\cite{CalinChangMarkina1, ChangMarkinaVasiliev}.

\begin{proposition}
The horizontality property is invariant under the left translation.
\end{proposition}

\begin{proof}
It can be shown that~\eqref{coord} does not change under the left translation. This implies the conclusion of the
proposition.
\end{proof}

\section{Hamiltonian system}\label{Sec3}

Once we have a system of curves, in our case the system of horizontal curves, we can define their length as in the Riemannian geometry. Let $\gamma:[0,t]\to \mathbb S^3$ be a horizontal curve such that $\gamma(0)=x$,
$\gamma(t)=y$, then the length $l(\gamma)$ of $\gamma$ is defined as follows
\begin{equation}\label{length}
l(\gamma)=\int_0^t\langle\dot\gamma, \dot\gamma\rangle^{1/2}\,ds=\int_0^t\big(a^2(s)+b^2(s)\big)^{1/2}\,ds.
\end{equation} Now we are able to define the distance between two points $x$ and $y$ by minimizing 
integral~\eqref{length} or the corresponding energy integral
$\int_0^t\big(a^2(s)+b^2(s)\big)\,ds$ under the non-holonomic
constraint~\eqref{hc}. This is the Lagrangian approach. The Lagrangian
formalism was applied to study  sub-Riemannian geometry on
$\mathbb S^3$ in~~\cite{CalinChangMarkina1, HurtadoRosales}. In  Riemannian geometry the
minimizing curve locally coinsides with the geodesic, but it is not
the case for  sub-Riemannian manifolds. Interesting examples
and discussions can be found, for instance,
in~\cite{LiuSussman,Montgomery1,Montgomery2,Montgomery,Strichartz}. Given the sub-Riemannian metric we
can form the Hamiltonian function defined on the cotangent bundle of
$\mathbb S^3$. A geodesic on a sub-Riemannian manifold is
defined as the  projection of a solution to the
corresponding Hamiltonian system onto the manifold. It is a good generalization of
the Riemannian case in the following sense. The Riemannian
geodesics (which are defined as curves with vanishing acceleration) can be lifted
to the solutions of the Hamiltonian system on the cotangent bundle.

Let us construct and describe sub-Riemannian geodesics
on $(\mathbb S^3,\mathcal D,\langle\cdot,\cdot\rangle)$. The left-invariant vector fields $X,Y,Z$ can be written using the matrices
\begin{equation*}
\label{eq:matrix}
I_1=\left[\array{rrrr}0 & 0 & -1 & 0
\\ 0 & 0 & 0 & -1
\\
1 & 0 & 0 & 0
\\ 0 & 1 & 0 & 0\endarray\right],\quad
I_2=\left[\array{rrrr}0 & 0 & 0 & -1
\\ 0 & 0 & 1 & 0
\\
0 & -1 & 0 & 0
\\ 1 & 0 & 0 & 0\endarray\right],\quad
I_3=\left[\array{rrrr}0 & -1 & 0 & 0
\\ 1 & 0 & 0 & 0
\\
0 & 0 & 0 & 1
\\ 0 & 0 & -1 & 0\endarray\right].
\end{equation*} In fact, $$X=\langle I_1x,\nabla x\rangle,\quad Y=\langle I_2x,\nabla x\rangle,
\quad Z=\langle I_3x,\nabla x\rangle.$$
The Hamiltonian  is defined as
$$H=\frac{1}{2}(X^2+Y^2)=\frac{1}{2}\Big(\langle I_1x,\xi\rangle^2
+\langle I_2x,\xi\rangle^2\Big),$$
or
\begin{equation}\label{Ham}
H=\frac{1}{2}(-x_3\xi_1-x_4\xi_2+x_1\xi_3+x_2\xi_4)^2+\frac{1}{2}(-x_4\xi_1+x_3\xi_2-x_2\xi_3+x_1\xi_4)^2,
\end{equation}
 where $\xi=\nabla x$. Then the Hamiltonian system follows as
\begin{equation}\label{hs}
\begin{split}
&\dot x=\frac{\partial H}{\partial \xi}\ \ \Rightarrow \ \ \dot
x=\langle I_1x,\xi\rangle\cdot (I_1x)+\langle I_2x,\xi\rangle\cdot (I_2x)\\
&\dot \xi=-\frac{\partial H}{\partial x}\ \ \Rightarrow \ \
\dot\xi=\langle I_1x,\xi\rangle\cdot (I_1\xi)+\langle
I_2x,\xi\rangle\cdot (I_2\xi).
\end{split}\end{equation}
As it was mentioned, a geodesic is the projection of a solution to the
Hamiltonian system onto the $x$-space. We obtain the following
properties.
\begin{itemize}
 \item [1.]{
Since $\langle I_1x,x\rangle=\langle I_2x,
x\rangle=\langle I_3x,x\rangle=0$, multiplying the first equation
of~\eqref{hs} by $x$, we get
$$\langle\dot x,x\rangle=0\ \ \Rightarrow \ \ |x|^2=const.$$
This asserts that {\it any solution to the Hamiltonian system belongs to the
sphere.} Taking the constant equal to $1$ we get geodesics on
$\mathbb S^3$.}
\item[2.]{Multiplying the first equation of~\eqref{hs} by
$I_3x$, we get
\begin{equation}\label{hc1}
\langle\dot x,I_3x\rangle=0,
\end{equation}
by the rule of multiplication for $I_1$, $I_2$, and $I_3$. The
reader easily recognizes the horizontality condition $\langle\dot
x,Z\rangle=0$ in~\eqref{hc1}. It means that {\it any solution to
the Hamiltonian system is a horizontal curve}}.
\item[3.]{Multiplying the first equation of~\eqref{hs} by $I_1x$, and then by $I_2x$,
we get
$$\langle \xi,I_1x\rangle=\langle\dot x,I_1x\rangle,
\qquad \langle \xi,xI_2\rangle=\langle\dot x,xI_2\rangle.$$ On the other hand, we know that $\langle\dot x,I_1x\rangle=a$ and
$\langle\dot x,xI_2\rangle=b$. The Hamiltonian  can be
written in the form
$$H=\frac{1}{2}\Big(\langle I_1x,\xi\rangle^2+\langle I_2x,\xi\rangle^2\Big)
=\frac{1}{2}\Big(\langle I_1x,\dot x\rangle^2+\langle I_2x,\dot
x\rangle^2\Big)=\frac{1}{2}\Big(a^2+b^2\Big).
$$
Thus, {\it the Hamiltonian  gives the kinetic energy $H=\frac{|\dot
q|^2}{2}$ which is  constant along the geodesics.}}
\item[4.]
{If we multiply the first equation of~\eqref{hs}
by $\dot x$, then we get
\begin{equation*}
|\dot x|^2=\langle I_1x,\xi\rangle^2+\langle I_2x,\xi\rangle^2
=\langle I_1x,\dot x\rangle^2+\langle I_2x,\dot
x\rangle^2=a^2+b^2=2H.
\end{equation*}  Therefore \begin{equation}\label{eq1}|\dot
x|^2=a^2+b^2.\end{equation}}
\end{itemize}

The following theorem was proved in \cite{ChangMarkinaVasiliev} and \cite{HurtadoRosales}. 

\begin{theorem}\label{th1}
The set of geodesics with constant velocity coordinates starting from the point $(1,0,0,0)$ forms the
unit sphere $\mathbb S^2$ in $\mathbb R^4$ parametrized as 
\[
(\cos s, \,\,0,\,\, \cos \psi \,\sin s, \,\,\sin\psi \, \sin s), \quad s\in [0,\pi], \quad \psi\in [0,2\pi).
\]
The integral line corresponding to the vertical vector field $Z$ starting from the point $(1,0,0,0)$ is parametrized as $(\cos\omega, \sin\omega,0,0)$, $\omega\in [0,2\pi)$.
\end{theorem}

\begin{remark}
For the arbitrary reference point the horizontal geodesics are parametrized by
\[
x(s)=x_0\cos s+(I_1\cos\psi +I_2 \sin\psi)x_0\sin s,
\]
and the vertical line by
\[
x(s)=x_0\cos s+I_3x_0\sin s,
\]
see \cite{ChangMarkinaVasiliev}.
\end{remark}

\section{Optimal control viewpoint}

The above Hamiltonian system and calculation of geodesics admits the optimal control interpretation.
The interplay of the control theory and sub-Riemannian geometry has been well known since early 80s. One of the pioneering contributions was made by Brockett \cite{Brockett}. He considered a time optimal control problem leading to the sub-Riemannian geometry in~$\mathbb R^3$, or to the Heisenberg group.
His results then were generalized in several ways, see e.g.,~\cite{Monroy}. Several results, already known by this time due to the fundamental Gaveau's work~\cite{Gaveau}, were rediscovered and 
the problem of finding normal and abnormal geodesics was formulated in terms of the optimal control, see e.g., \cite{Agrachev2, LiuSussman}.
Pontryagin's maximum principle provides such optimal controls. Interesting features of such Hamiltonian systems are symmetries given by the first integrals although such systems generally are not (Frobenius) integrable
because of singular geometric background, i.e., constraints on the velocities can not be re-written in terms of the configuration coordinates. A good reference to the control theory viewpoint is~\cite{Bloch}.

Let us consider  the following time optimal control problem given by the system
\begin{equation}\label{pseudo1}
\begin{array}{lll}
\dot{x}_1 & = & -ux_3-v x_4,   \\
\dot{x}_2 & = & -ux_4+v x_3,   \\
\dot{x}_3 & = & ux_1-v x_2,   \\
\dot{x}_4 & = & ux_2+v x_1,   
\end{array}
\end{equation}
with the cost functional
\[
E=\frac{1}{2}\int_0^t\langle \bu, \bu \rangle ds,
\]
where $\bu= (u,v)$. The functional $E$ represents the total kinetic energy. The system is encoded in the kernel of the contact 1-form \eqref{forma}.

The pseudo-Hamiltonian given by the Pontryagin Maximum Principle  for this system admits the form
\begin{equation}\label{pseudo2}
\mathcal H=-\frac{1}{2}(u^2+v^2)+u(-x_3\xi_1-x_4\xi_2+x_1\xi_3+x_2\xi_4)+v(-x_4\xi_1+x_3\xi_2-x_2\xi_3+x_1\xi_4),
\end{equation}
and the system for covectors becomes
\begin{equation}\label{pseudo3}
\begin{array}{lll}
\dot{\xi}_1 & = & -u\xi_3-v \xi_4,   \\
\dot{\xi}_2 & = & -u\xi_4+v \xi_3,   \\
\dot{\xi}_3 & = & u\xi_1-v \xi_2,   \\
\dot{\xi}_4 & = & u\xi_2+v \xi_1.   
\end{array}
\end{equation}
The system (\ref{pseudo1}--\ref{pseudo3}) for position coordinates may be rewritten in the following form
\begin{equation*}
\begin{array}{lll}
u & = & -x_3\dot{x}_1-x_4\dot{x}_2+x_1\dot{x}_3+x_2\dot{x}_4,   \\
v & = & -x_4\dot{x}_1+x_3\dot{x}_2-x_2\dot{x}_3+x_1\dot{x}_4,   \\
0 & = & -x_2\dot{x}_1+x_1\dot{x}_2+x_4\dot{x}_3-x_3\dot{x}_4 ,  \\
0 & = & x_1\dot{x}_1+x_2\dot{x}_2+x_3\dot{x}_3+x_4\dot{x}_4,   
\end{array}
\end{equation*}
which has a clear geometric meaning. Indeed, $u$ and $v$ are the coefficients of the velocity vector
$uX+vY$, the third equation is just the horizontality condition and the fourth means that the trajectory
belongs to a sphere. 

From the Hamiltonian system one derives four first integrals
\begin{equation*}
\begin{array}{lll}
J_1 & = & x_1\xi_1+x_2\xi_2+x_3\xi_3+x_4\xi_4,   \\
J_2 & = & -x_2\xi_1+x_1\xi_2+x_4\xi_3-x_3\xi_4,   \\
J_3 & = & -x_3\xi_1+x_4\xi_2+x_1\xi_3-x_2\xi_4,  \\
J_4 & = & -x_4\xi_1-x_3\xi_2+x_2\xi_3+x_1\xi_4.    
\end{array}
\end{equation*}
The Poisson structure is given by the Poisson brackets
\[
[F,G]=\sum\limits_{k=1}^4 \frac{\partial F}{\partial x_k}\frac{\partial G}{\partial \xi_k}- \frac{\partial G}{\partial x_k}\frac{\partial F}{\partial \xi_k}.
\] 
The integrals $J_1$ and $J_2$ represent natural symmetries (following two natural geometric conditions: $J_1$ is the normal covector and $J_2$ gives the horizontality condition) and $J_3$, $J_4$ give hidden symmetries.
All first integrals are involutive in pairs $[J_k, J_m]=0$, $k,m=1,\dots, 4$, which implies Liouville integrabilty of the above Hamiltonain system.  Observe, that the Hamiltonian system for the Heisenberg group is not Liouville integrable as well as the Hamiltonian system corresponding to sub-Riemannian geometry on $SO(n)$ for $n\geq 4$, see \cite{Boscain, Respondek}. Let us remark that the optimal control problem in the sub-Riemannian geometry on $SO(n)$ can be viewed as the problem of optimal laser-induced population transfer in $n$-level quantum systems, see \cite{Boscain}.

In order to find geodesics we can use the Pontryagin Maximum Principle \cite{Pontryagin} which states that any normal geodesic is a projection of a bicharacteristic which is a solution to the above Hamiltonian system on the cotangent bundle with the control $\bu^*$ which maximizes the pseudo-Hamiltonian $\mathcal H$, i.e., satisfies the equation 
\[
\frac{\partial \mathcal H}{\partial u}=\frac{\partial \mathcal H}{\partial v}=0.
\]
This problem is equivalent to the geometric problem of minimizing the Carnot-Carath\'e\-odory distance
(or, equivalently, sub-Riemannian energy) in the optimal control problem for our control-linear system. The optimal control admits the form
\[
u^*=-x_3\xi_1-x_4\xi_2+x_1\xi_3+x_2\xi_4, \quad v^*=-x_4\xi_1+x_3\xi_2-x_2\xi_3+x_1\xi_4.
\]
Substituting $\bu^*$ in the Hamiltonian system we obtain the geodesic equation \eqref{Ham} and \eqref{hs}. The importance of integrability of the sub-Riemannian geodesic equation was argued by Brockett and Dai \cite{BrocketDai}, who showed the explicit integrability in some special cases in terms of elliptic functions and discussed applications to controllability problems. But the question of integrabilty of Hamiltonain systems associated with nonholonomic distributions has a long history, see the survey \cite{Vershik}
for the historical account.

As it was shown in \cite{Montgomery1}, abnormal geodesics are not geometrically relevant for step 2 groups. Nevertheless, we give here independent treatment of abnormal geodesics from the Pontryagin Maximum Principle viewpoint. 
The pseudo-Hamiltonian in this case becomes
\[
\mathcal H_0=u(-x_3\xi_1-x_4\xi_2+x_1\xi_3+x_2\xi_4)+v(-x_4\xi_1+x_3\xi_2-x_2\xi_3+x_1\xi_4)=uJ_3+vJ_4.
\]
The Pontryagin Maximum Principle implies that $\mathcal H_0$ vanishes along the extremal. We can assume that the velocity coordinates $u$ and $v$ do not vanish simultaneously.
After differentiating $J_3$ and $J_4$ along the extremal we obtain
\[
0=\dot{J_3}=[J_3, \mathcal H_0]=-2vJ_2,
\]
\[
0=\dot{J_4}=[J_4, \mathcal H_0]=2uJ_2.
\]
Let us suppose that $u$ does not vanish on some time interval $s\in U$. Then, $J_2=0$ on this interval, and being the first integral, it is vanishing everywhere. Then we obtain
\[
0=\dot{J_2}=[J_2, \mathcal H_0]=-2uJ_4,
\]
and $J_4$ is identically 0 by the same reason. Therefore, $J_3\equiv 0$. Solving the system $J_k=0$, $k=1,\dots,4$, with respect to $x_k$ we see that the discriminant of this system is 1. Fixing initial conditions for the Hamiltonian system (\ref{pseudo1}--\ref{pseudo3}) we deduce that $\xi_k\equiv 0$, $k=1,\dots,4$, and only stationary solution is valid.

For normal geodesics we have that along the extremal
\[
\mathcal H=H=\frac{1}{2}(-x_3\xi_1-x_4\xi_2+x_1\xi_3+x_2\xi_4)^2+\frac{1}{2}(-x_4\xi_1+x_3\xi_2-x_2\xi_3+x_1\xi_4)^2,
\]
and is given as in \eqref{Ham}.

\section{Geodesics and modified action}\label{Modif.function}

\subsection{Cartesian coordinates}

Fix the initial point $x^{(0)}=(1,0,0,0)$. It is convenient to introduce complex coordinates $z=x_1+ix_2$, $w=x_3+ix_4$, $\varphi=\xi_1+i\xi_2$, and $\psi=\xi_3+i\xi_4$. Hence, the Hamiltonian $H$ admits the form $H=\frac{1}{2}|\bar w\varphi -z\bar\psi|^2$ (compare with \eqref{Ham}). The corresponding Hamiltonian system becomes
\begin{equation*}
\begin{array}{lllll}
\dot z & = &  w(\bar{w} \varphi- z\bar{\psi}),\quad & &z(0)=1, \\
\dot w & = & -z(w\bar{\varphi}-\bar{z}\psi), \quad & &w(0)=0,\\
\dot {\bar{\varphi}} & = & \bar{\psi}(w\bar{\varphi}-\bar{z}\psi), \quad & &\bar{\varphi}(0)=A-iB,\\
\dot {\bar{\psi}} & = & -\bar{\varphi}(\bar{w} \varphi-z\bar{\psi}), \quad & &\bar{\psi}(0)=C-iD,
\end{array}
\end{equation*}
and $H=\frac{1}{2}\re(\dot{z}\bar{\varphi}+\dot{w}\bar{\psi})$.
Here the constants $B,C$, and $D$ have the following dynamical meaning: $\dot w(0)=C+iD$,
and $B=-i \ddot w(0)/2\dot w(0)$ or if we write in real variables, $C=\dot x_3(0)$, $D=\dot x_4(0)$,
$B=\frac{1}{2}(\dot x_3(0)\ddot x_4(0)-\dot x_4(0)\ddot x_3(0))/(\dot{x}_3^2(0)+\dot{x}_4^2(0))$. If we denote 
\[
k=\frac{B}{\sqrt{C^2+D^2}},
\]
then $|k|$ is the curvature of a geodesic at the initial point.
This complex Hamiltonian system has the
first integrals
\[
\begin{array}{rcl}
z\psi -w \varphi  & = &  C+iD, \\
z\bar{\varphi}+w \bar{\psi} & = & A-iB,\\
\end{array}
\]
and we have $|z|^2+|w|^2=1$  as a
normalization. Therefore,
\[
\begin{array}{rcl}
\varphi  & = & z(A+iB)-\bar w ( C+iD), \\
\psi  & = & \bar z (C+iD) +w (A+iB).\\
\end{array}
\]
Let us introduce an auxiliary function $p=\bar w/{z}$. Then substituting $\varphi$ and $\psi$ in the
Hamiltonian system we get the equation for $p$ as
\[
\dot{p}=(C+iD)p^2-2iBp+(C-iD),\quad p(0)=0.
\]
The solution is
\[
p(s)=\frac{(C-iD)\sin(s\sqrt{B^2+C^2+D^2})}{\sqrt{B^2+C^2+D^2}\cos(s\sqrt{B^2+C^2+D^2})+iB \sin(s\sqrt{B^2+C^2+D^2})}.
\]
Taking into account that $\dot z \bar z= -w \dot{\bar w}$, we get the solution
\begin{equation}\label{geodz}
z(s)=\left(\cos(s\sqrt{B^2+C^2+D^2})+i\frac{B}{\sqrt{B^2+C^2+D^2}}\sin(s\sqrt{B^2+C^2+D^2})\right)e^{-iBs},
\end{equation}
and
\begin{equation}\label{geodw}
w(s)=\frac{C+iD}{\sqrt{B^2+C^2+D^2}} \sin(s\sqrt{B^2+C^2+D^2})e^{iBs}.
\end{equation}

\begin{remark}
Let us consider three limiting cases.
If $B=0$, then we get the solutions with constant horizontal velocity coordinates $$z(s)=\cos s,
\qquad w(s)=(\dot x_3(0)+i\dot x_4(0))\sin s$$ which lie on the horizontal 2-sphere, and  a geodesic joining two given points on it is unique.
If $C^2+D^2=0$, then the only solution $w(s)$ to the Hamiltonian system is $w(s)\equiv 0$. The horizontality condition in this case is read as $x_2\dot{x}_1=x_1\dot{x}_2$, and the solution is a straight line which contradicts the condition $|z|^2=1$. So $H=\frac{1}{2}(C^2+D^2)>0$.
\end{remark}
 
Now we want to find geodesics joining two given points.

\begin{theorem}\label{numbergeod1}
Let $Q$ be a point of the vertical line, i.~e. $Q=(\cos \omega,\sin \omega, 0,0)$, $\omega\in(-\pi,0)\cup(0,\pi)$, then there are countably
many geometrically different geodesics $\gamma_n$ connecting $P=(1,0,0,0)$ with $Q$. They have the following parametric representation
\begin{eqnarray}\label{geod1}
 z_n(s) & = & \Big(\cos(s\frac{\pi n}{t})-
 i\frac{\omega}{\pi n}\sin(s\frac{\pi n}{t})\Big)e^{-\frac{is\omega}{t}},\\
w_n(s) & = & (C+iD)\frac{t}{\pi n}
\sin(s\frac{\pi n}{t})e^{\frac{is\omega}{t}},\nonumber
\end{eqnarray}
$n\in\mathbb Z\setminus\{0\}$, $s\in[0,t]$, and the length of geodesics 
$\gamma_n$ is given as $l_n=\frac{1}{\sqrt{2}}\sqrt{(\pi n)^2-\omega^2}$.
\end{theorem}
\begin{proof}
The geodesics are parametrized in the time interval $s\in [0,t]$.
 If the point $Q=(z(t),w(t))=(z,w)$ belongs to the vertical line starting at
$P=(1,0,0,0)$, then $|z|=1$ and $|w|=0$ provided that
$-Bt=\omega$, in what follows,  \[\cos^2
(t\sqrt{B^2+C^2+D^2})+\frac{B^2}{B^2+C^2+D^2}\sin^2(t\sqrt{B^2+C^2+D^2})=1,
\]
\[
\sin(t\sqrt{B^2+C^2+D^2})=0,\qquad -Bt=\omega.\] 
These equations imply
\begin{equation}\label{relations}
t=\frac{\pi n}{\sqrt{B^2+C^2+D^2}}>0,\quad -Bt=\omega.
\end{equation}
The latter relations give
\[
B^2=B^2_n\equiv \frac{\omega^2(C^2+D^2)}{(\pi n)^2-\omega^2}.
\]
Substituting \eqref{relations} in the solutions to the Hamiltonian system we come to the parametric
representation given in the formulation of the theorem. 
The first relation of \eqref{relations} yields
\[
\sqrt{C^2+D^2}=\frac{1}{t}\sqrt{(\pi n)^2-\omega^2}.
\]
The length of each geodesic is given as
\[
l_n=t\sqrt{H}=\frac{t}{\sqrt{2}}\sqrt{C^2+D^2}=\frac{1}{\sqrt{2}}\sqrt{(\pi n)^2-\omega^2}.
\]
This finishes the proof.
\end{proof}

\begin{remark}
In the formulation of the theorem the words `geometrically different' mean that due to
the rotation of the argument of $C+iD$ in $w(s)$, there exist uncountably many geodesics.
\end{remark}

So far we have had a clear picture of trivial geodesics whose velocity
has constant coordinates. They are essentially unique (up to
periodicity). The situation with geodesics joining the point
$(1,0,0,0)$ with the points of the vertical line has been
described in the preceding theorem. Let us consider the generic
position of the right endpoint $(z,w)$, $z=re^{i\xi_1}$, $w=\rho e^{i\xi_2}$ on $\mathbb S^3$.

\begin{remark}
First we consider three limiting cases. If $\rho=1$, then  $B=0$ and $r=0$, and the point
lies on the horizontal 2-sphere. If $\rho=0$, then $\sin(t\sqrt{B^2+C^2+D^2})=0$ and the point $z=\pm(\cos (Bt)-i\sin(Bt))$, $w=0$ belongs to the vertical line.  If $\arg z=0,\pi$, then $z=\pm r$, $w=\pm \sqrt{1-r^2}(\cos \xi_2+i\sin\xi_2)$ is
 a point on the horizontal 2-sphere.
\end{remark}

In other situations we have the following theorem.

\begin{theorem}\label{numbergeod2}
Given an arbitrary point $(z,w)\in \mathbb S^3$ which neither belongs to the vertical line  nor to
the horizontal sphere $\mathbb S^2$, there is a finite number of geometrically different geodesics
joining the initial point (north  pole) $P=(1,0,0,0)\in \mathbb S^3$ with $Q=(z,w)$.
\end{theorem}
\begin{proof}
Let us denote
\[
z=re^{i\xi_1},\quad w=\rho e^{i\xi_2},\quad C+iD=\sqrt{C^2+D^2}e^{i\theta}.
\]
Then from \eqref{geodz} and \eqref{geodw} we have that
\begin{equation}\label{ugol}
\rho^2=\frac{C^2+D^2}{B^2+C^2+D^2}\sin^2(t\sqrt{B^2+C^2+D^2}),\quad \mbox{and\  \ }\xi_2=Bt+\theta,
\end{equation} where $t$ is the right end of the time interval $s\in[0,t]$ at which the endpoint $Q$ is reached. 
We suppose for the moment that the angles $s\sqrt{B^2+C^2+D^2}$ and $tB$ are from the first quadrant.
Other cases are treated similarly. Then we have
\[
z=\left(\sqrt{1-\frac{B^2+C^2+D^2}{C^2+D^2} \rho^2}+i\frac{B\rho}{\sqrt{C^2+D^2}}\right)e^{i(\theta-\xi_2)},
\]
and
\[
\xi_1=\theta-\xi_2+\arctan\frac{B\rho}{\sqrt{C^2+D^2-(B^2+C^2+D^2)\rho^2}}.
\]
The first expression in \eqref{ugol} leads to the value of the length parameter $t$ as
\[
t=\frac{1}{\sqrt{B^2+C^2+D^2}}\arcsin \left(\rho\sqrt{1+\frac{B^2}{C^2+D^2}}\right),
\]
and the second to
\[
\xi_2=\theta+\frac{B}{\sqrt{B^2+C^2+D^2}}\arcsin \left(\rho\sqrt{1+\frac{B^2}{C^2+D^2}}\right).
\]
Substituting $\theta$ in the latter equation we come to an equation which depends only on 
\[
k=\frac{B}{\sqrt{C^2+D^2}},
\]
which we rewrite as
\begin{equation}\label{equa}
\sin\left(\sqrt{1+\frac{1}{k^2}}\Bigg[\arctan\Big(\frac{k\rho}{1-(1+k^2)\rho^2}\Big)-\xi_1\Bigg]\right)=\rho\sqrt{1+k^2},
\end{equation}
or as an equation for the parameter $k$, which is the curvature of the geodesic at the initial moment. We through away the trivial cases $k=0$ and $\xi_1=0$ excluded from the theorem (see the remark before the theorem).  

Observe that $\theta-\xi_2=\xi_1-\arctan\Big(\frac{k\rho}{1-(1+k^2)\rho^2}\Big)$ is non-vanishing because $B\neq 0$ from \eqref{ugol}.  So the left-hand side of equation \eqref{equa} is a  function of $k$ which is bounded by 1 in
absolute value and fast oscillating about the point $k=0$. The right-hand side of \eqref{equa}
is an even function increasing for $k>0$, see Figure~\ref{fig2}. Therefore, there exists a countable number of non-vanishing
different solutions $\{k_n\}$ of the equation \eqref{equa} within the interval $|k|\leq \sqrt{\frac{1}{\rho^2}-1}=\frac{|z|}{|w|}$
with a limit point at the origin. 
\begin{figure}[ht]
\centering \scalebox{0.3}{\includegraphics{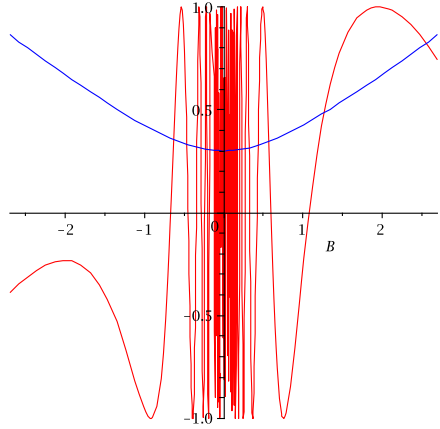}}
\caption[]{Solutions to the equation~\eqref{equa}  }\label{fig2}
\end{figure}

However, in order to define all parameters $B$, $C$, and $D$ we need
to solve the equations (\ref{ugol}), (\ref{equa}), and not all $k_n$ satisfy all three equations. Let us consider positive $k_n$.
We calculate the argument of $z$ as
\[
\xi_1=-Bt+\arctan\Bigg[\frac{B}{\sqrt{B^2+C^2+D^2}}\tan\left(t\sqrt{B^2+C^2+D^2}\right)\Bigg]
\]
\[
=-Bt+\arctan\Bigg[\frac{k_n\rho}{\sqrt{1-(1+k_n^2)\rho^2}}\Bigg]
\]
\[
<-k_n t \sqrt{C^2+D^2}+\frac{k_n\rho}{\sqrt{1-(1+k_n^2)\rho^2}}.
\]

On the other hand, we have
\[
\sqrt{C^2+D^2}=\frac{\arcsin(\rho\sqrt{1+k_n^2})}{t\sqrt{1+k_n^2}}>\frac{\rho}{t}.
\]
Observe that due to the remark before this theorem, $\xi_1>0$ and $0<\rho<1$.
Therefore, we deduce the inequality
\[
\xi_1<k_n\rho\frac{1-\sqrt{1-\rho^2(1+k_n^2)}}{\sqrt{1-\rho^2(1+k_n^2)}},
\]
or
\begin{equation}\label{equ}
k_n\rho>\xi_1\frac{\sqrt{1-\rho^2(1+k_n^2)}}{1-\sqrt{1-\rho^2(1+k_n^2)}}.
\end{equation}
The right-hand side of the inequality \eqref{equ} decreases with respect to $k_n>0$.

Set $\varepsilon=\frac{1+\rho^2}{2}$. If $\varepsilon <\rho^2(1+k_n^2)<1$, then immediately
we have the inequality
$k_n^2>\frac{1}{2}(\frac{1}{\rho^2}-1)>0$. If $0 <\rho^2(1+k_n^2)\leq \varepsilon$, then the inequality
\eqref{equ} implies that
\[
k_n>\xi_1\frac{\sqrt{1-\varepsilon}}{\rho(1-\sqrt{1-\varepsilon})}=\xi_1\frac{\sqrt{1-\rho^2}}{\rho(\sqrt{2}-\sqrt{1-\rho^2})}>0.
\]
Finally, we obtain
\[
k_n>\min\left\{\xi_1\frac{\sqrt{1-\rho^2}}{\rho(\sqrt{2}-\sqrt{1-\rho^2})}, \sqrt{\frac{1}{2}(\frac{1}{\rho^2}-1)}\right\}\equiv b(\xi_1,\rho)>0.
\]
This proves that all positive solutions to the equation \eqref{equa} must belong to the interval $(b(\xi_1,\rho), \sqrt{\frac{1}{\rho^2}-1})$, hence there are only finite number of such $k_n$. The same arguments are applied for negative values of $k_n$. 
\end{proof}

\begin{remark}
If $\rho$ is approaching $0$, the point $Q$ is approaching the vertical line and the value of $k_n$ becomes
\[
k_n=\frac{\pm\xi_1}{\sqrt{(\pi n)^2-\xi_1^2}},
\]
and the solution is reduced to the case considered in Theorem \ref{numbergeod1} with $\omega=\xi_1$, i.e., the number of geodesics is increasing infintely.
\end{remark}
\begin{remark}
Given two  points $P$ and $Q$, we find the initial velocity by equation~\eqref{ugol} and the initial curvature by equation~\eqref{equa}.
\end{remark}

\subsection{Hyperspherical coordinates}
 
Let us use now the hyperspherical coordinates 
\begin{eqnarray}\label{hspc}
x_1+ix_2 & = & e^{i\zeta_1}\cos\eta,\\
x_3+ix_4 & = & e^{i\zeta_2}\sin\eta,\qquad \eta\in(0,\pi/2),\quad\zeta_1,\zeta_2\in[-\pi,\pi),\nonumber
\end{eqnarray}
to write the Hamiltonian system.

The horizontal coordinates are written as
\begin{eqnarray*}
a & = & \dot\eta\cos(\zeta_1-\zeta_2)+(\dot\zeta_1+\dot\zeta_2)\sin(\zeta_1-\zeta_2)\frac{\sin 2\eta}{2},\\
b & = & -\dot\eta\sin(\zeta_1-\zeta_2)+(\dot\zeta_1+\dot\zeta_2)\cos(\zeta_1-\zeta_2)\frac{\sin 2\eta}{2},\\
c & = & \dot\zeta_1\cos^2\eta-\dot\zeta_2\sin^2\eta.
\end{eqnarray*}

The horizontality condition in hyperspherical coordinates becomes
$$\dot\zeta_1\cos^2\eta-\dot\zeta_2\sin^2\eta=0.$$

The horizontal 2-sphere in Theorem~\ref{th1} is obtained from the
parametrization~\eqref{hspc}, if we set $\zeta_1=0$, $\zeta_2=\psi$,
$\eta=s$ or $\eta=\pi-s$.  The vertical line is obtained from the
parametrization~\eqref{hspc} setting $\eta=0$, $\zeta_1=s$.

Writing the vector fields $N,Z,X,Y$ in the hyperspherical
coordinates we get $$N=-2\cotan2\eta\partial_{\eta},\quad
Z=\partial_{\zeta_1}-\partial_{\zeta_2},$$
$$X=\sin(\zeta_1-\zeta_2)\tan\eta\partial_{\zeta_1}+\sin(\zeta_1-\zeta_2)\cotan\eta\partial_{\zeta_2}
+2\cos(\zeta_1-\zeta_2)\partial_{\eta},$$
$$Y=\cos(\zeta_1-\zeta_2)\tan\eta\partial_{\zeta_1}
+\cos(\zeta_1-\zeta_2)\cotan\eta\partial_{\zeta_2}-2\sin(\zeta_1-\zeta_2)\partial_{\eta}.$$
In this parametrization some similarity with the Heisenberg group
can be shown. The commutator of two horizontal vector fields $X,Y$ gives the constant vector field $Z$ which is orthogonal to the
horizontal vector fields at each point of the manifold. In
hyperspherical coordinates it is easy to see that the form
$\omega=\cos^2\eta d\zeta_1-\sin^2\eta d\zeta_2$, that defines the
horizontal distribution is contact because 
$$\omega\wedge d\omega=\sin (2\eta)\,d\eta\wedge d\zeta_1\wedge
d\zeta_2=2dV,$$ where $dV$ is the volume form. The sub-Laplacian is
defined as
$$\frac{1}{2}(X^2+Y^2)=\frac{1}{2}(\tan^2\eta\partial^2_{\zeta_1}+\cotan^2\eta\partial^2_{\zeta_2}
+4\partial^2_{\eta}+2\partial_{\zeta_1}\partial_{\zeta_2}).$$
The principal symbol is given by the Hamiltonian
$$H(\zeta_1,\zeta_2,\eta,\psi_1,\psi_2,\theta)=\frac{1}{2}(\tan^2\eta\psi^2_1+\cotan^2\eta\psi^2_2
+4\theta^2+2\psi_1\psi_2),$$
with the covectors $\psi_k\sim\partial_{\zeta_k}$, $k=1,2$, $\theta=\partial_{\eta}$.
It gives the Hamiltonian system 
\begin{equation}\label{Hsystem}
\begin{array}{lll}
\medskip\dot\zeta_1 & = & \frac{\partial H}{\partial \psi_1}=\psi_1\tan^2\eta+\psi_2\\
\medskip\dot\zeta_2 & = & \frac{\partial H}{\partial \psi_2}=\psi_2\cotan^2\eta+\psi_1\\
\medskip\dot\eta & = & \frac{\partial H}{\partial \theta}=4\theta\\
\medskip\dot\psi_1 & = & -\frac{\partial H}{\partial \zeta_1}=0\\
\medskip\dot\psi_2 & = & -\frac{\partial H}{\partial \zeta_2}=0\\
\medskip\dot\theta & = & -\frac{\partial H}{\partial \eta}=-\psi_1^2\frac{\tan\eta}{\cos^2\eta}
+\psi_2^2\frac{\cotan\eta}{\sin^2\eta}.
	\end{array}
\end{equation} 

\subsection{Geodesics in hyperspherical coordinates}
Let us find geodesics \[\gamma(s)=(\zeta_1(s),\zeta_2(s),\eta(s)), \quad s\in [0,t]\] joining the points $P=\gamma(0)=(\zeta_1^0,\zeta_2^0,\eta_0)$ and $Q=\gamma(t)=(\zeta_1,\zeta_2,\eta)$. They are obtained as projections  of the solutions to system \eqref {Hsystem} onto the sphere.

Observe that the system \eqref{Hsystem} is coupled and the system
\begin{equation}\label{Hsystem2}
\begin{array}{lll}
\dot\eta & = & \frac{\partial H}{\partial \theta}=4\theta\\
\dot\theta & = & -\frac{\partial H}{\partial \eta}=-\psi_1^2\frac{\tan\eta}{\cos^2\eta}
+\psi_2^2\frac{\cotan\eta}{\sin^2\eta}
\end{array}
\end{equation} 
with the boundary conditions $\eta(0)=\eta_0$,  $\eta(t)=\eta$,  is independent.

Multiplying the equations of this system crosswise we obtain
\[
4\theta\dot{\theta}=\left(-\psi_1^2\frac{\tan\eta}{\cos^2\eta}
+\psi_2^2\frac{\cotan\eta}{\sin^2\eta}\right)\dot{\eta},
\]
or
\[
\frac{d}{dt}(2\theta^2)=\frac{d}{dt}(-\frac{1}{2}\psi_1^2\tan^2\eta-\frac{1}{2}\psi_2^2\cot^2\eta).
\]
Therefore,
\begin{equation}\label{eqtheta}
4\theta^2(s)=4\theta^2_0-\psi_1^2\tan^2\eta(s)-\psi_2^2\cot^2\eta(s)+\psi_1^2\tan^2\eta_0+\psi_2^2\cot^2\eta_0.
\end{equation}
The constant $\theta_0$ is a constant of integration which will be further expressed in terms of boundary conditions for $P$ and $Q$.
Let us substitute the expression for $\theta(s)$ in the second equation of the system~\eqref{Hsystem2}. We obtain
\begin{equation}\label{diffeq}
\dot{\eta}=4\theta=\pm2\sqrt{4\theta^2_0+\psi_1^2\tan^2\eta_0+\psi_2^2\cot^2\eta_0-\psi_1^2\tan^2\eta(s)-\psi_2^2\cot^2\eta(s)}.
\end{equation}
Observe that the expression under the square root is non-negative for all $s\in [0,t]$. Let us consider the case of increasing $\eta$ and (+) in front of the square root (which is assumed to be positive). Negative case will be treated later. Changing variables $u=\sin^2\eta\in[0,1]$, we arrive at
\[
\frac{1}{4}\dot{u}=\sqrt{u(1-u)(4\theta_0^2+(\psi_1\tan\eta_0+\psi_2\cot\eta_0)^2)-(\psi_1 u+\psi_2(1-u))^2}.
\]
The square polynomial under the root is reduced to
\[
\left(\frac{\psi_1^2}{\cos^2\eta_0}+\frac{\psi_2^2}{\sin^2\eta_0}+4\theta_0^2\right)\left\{
\frac{\left(\frac{\psi_1^2}{\cos^2\eta_0}+\frac{\psi_2^2}{\sin^2\eta_0}+4\theta_0^2+\psi_2^2-\psi_1^2\right)^2}{4\left(\frac{\psi_1^2}{\cos^2\eta_0}+\frac{\psi_2^2}{\sin^2\eta_0}+4\theta_0^2\right)^2}
-\frac{\psi_2^2}{\frac{\psi_1^2}{\cos^2\eta_0}+\frac{\psi_2^2}{\sin^2\eta_0}+4\theta_0^2}-\right.
\] 
\[
\left. -\Bigg[u-\frac{\frac{\psi_1^2}{\cos^2\eta_0}+\frac{\psi_2^2}{\sin^2\eta_0}+4\theta_0^2+\psi_2^2-\psi_1^2}{2(\frac{\psi_1^2}{\cos^2\eta_0}+\frac{\psi_2^2}{\sin^2\eta_0}+4\theta_0^2)}\Bigg]^2\right\}.
\]
The polynomial is non-negative for all $u\in[0,1]$, as it was mentioned before. Therefore,
\[
\frac{\left(\frac{\psi_1^2}{\cos^2\eta_0}+\frac{\psi_2^2}{\sin^2\eta_0}+4\theta_0^2+\psi_2^2-\psi_1^2\right)^2}{4\left(\frac{\psi_1^2}{\cos^2\eta_0}+\frac{\psi_2^2}{\sin^2\eta_0}+4\theta_0^2\right)^2}
-\frac{\psi_2^2}{\frac{\psi_1^2}{\cos^2\eta_0}+\frac{\psi_2^2}{\sin^2\eta_0}+4\theta_0^2}
\] 
is non-negative too. Integrating \eqref{diffeq} gives us
\begin{equation}\label{eqeta}
\frac{\sin^2\eta(s)-\frac{1}{2}\left(1+\frac{\psi_2^2-\psi_1^2}{A}\right)}{\sqrt{\frac{1}{4}\left(1+\frac{\psi_2^2-\psi_1^2}{A}\right)^2-\frac{\psi_2^2}{A}}}=\sin(4s\sqrt{A}+const),
\end{equation}
where we introduce the notation
\begin{equation}\label{A}
A=\frac{\psi_1^2}{\cos^2\eta_0}+\frac{\psi_2^2}{\sin^2\eta_0}+4\theta_0^2>0.
\end{equation}
It is convenient to use the constants $\tilde\psi_1=\psi_1/\sqrt{A}$, $\tilde\psi_2=\psi_2/\sqrt{A}$, and $A$ instead of $\psi_1$, $\psi_2$, and $\theta_0$.
Now we consider the sign (-) in front of the square root. Finally, our solution is written as
\begin{equation}\label{eq:eta}
\frac{\sin^2\eta(s)-\frac{1}{2}(1+\tilde\psi_2^2-\tilde\psi_1^2)}{\sqrt{\frac{1}{4}(1+\tilde\psi_2^2-\tilde\psi_1^2)^2-\tilde\psi_2^2}}=\sin(\pm 4s\sqrt{A}+const),
\end{equation}
where
\[
const=\arcsin\frac{\sin^2\eta_0-\frac{1}{2}(1+\tilde\psi_2^2-\tilde\psi_1^2)}{\sqrt{\frac{1}{4}(1+\tilde\psi_2^2-\tilde\psi_1^2)^2-\tilde\psi_2^2}}+2\pi n.
\]

Let us turn to the solution of the boundary value problem with the boundary conditions $\zeta_1(0)=\zeta_2(0)=0$, $\eta(0)=\pi/4$, and $\zeta_1(t)=\zeta_1$, $\zeta_2(t)=\zeta_2$, $\eta(t)=\eta$.  Observe that the chosen parametrization does not give us a chart about the north pole $(1,0,0,0)$. So we can shift the considerations by a left-invariant group action to any initial point, e.g., $x_0=(1/\sqrt{2},0,1/\sqrt{2},0)$.  
The horizontal geodesics starting from the point $x_0$ admit the form
\begin{eqnarray*}
x_1(s)+ix_2(s)=\frac{1}{\sqrt{2}}\left(\cos s - \cos\sigma \sin s +i \sin\sigma\sin s\right),\\
x_3(s)+ix_4(s)=\frac{1}{\sqrt{2}}\left(\cos s + \cos\sigma \sin s +i \sin\sigma\sin s\right),
\end{eqnarray*}
see Theorem~\ref{th1} and the remark thereafter. In the hyperspherical coordinates it looks as
\begin{eqnarray*}
\cos\eta(s)=\sqrt{\frac{1}{2}-\cos\sigma\sin s\cos s},\\
\sin\zeta_1(s)=\frac{\sin\sigma\sin s}{\sqrt{1-\cos\sigma\sin 2s}},\\
\sin\zeta_2(s)=\frac{\sin\sigma\sin s}{\sqrt{1+\cos\sigma\sin 2s}},
\end{eqnarray*}
where $\sigma$ is some constant from the interval $[0,\pi]$. The horizontal surface is given by the relation $\sin\zeta_1=\sin\zeta_2\tan\eta$.
The vertical line is written as 
\[
\eta(s)\equiv \frac{\pi}{4},\quad \zeta_1=s,\quad \zeta_2=-s.
\]

Substituting $s=0$ and $s=t$ gives us the expression of $A=A^{(n)}(\eta_0,\eta,\tilde\psi_1,\tilde\psi_2,t)$ as
\begin{eqnarray}\label{eqA}
A &=& \frac{1}{16t^2}\left(\arcsin \frac{\sin^2\eta-\frac{1}{2}(1+\tilde\psi_2^2-\tilde\psi_1^2)}{\sqrt{\frac{1}{4}(1+\tilde\psi_2^2-\tilde\psi_1^2)^2-\tilde\psi_2^2}}
\right.\\
 &-& \left.\arcsin\frac{\sin^2\eta_0-\frac{1}{2}(1+\tilde\psi_2^2-\tilde\psi_1^2)}{\sqrt{\frac{1}{4}(1+\tilde\psi_2^2-\tilde\psi_1^2)^2-\tilde\psi_2^2}}+2\pi n \right)^2.\nonumber
\end{eqnarray}

Then we find the parametric representation for the functions $\zeta_1(s)$ and $\zeta_2(s)$. Let us proceed by writing
\[
\tan^2\eta(s)=-1+\frac{1}{1-\frac{1}{2}\left(1+\tilde\psi_2^2-\tilde\psi_1^2\right)- D_0\sin( \pm4s\sqrt{A}+D_1)},
\]
where $D_0$ and $D_1$ are the constants
\[
D_0=\sqrt{\frac{1}{4}\left(1+\tilde\psi_2^2-\tilde\psi_1^2\right)^2-\tilde\psi_2^2},
\]
\[
D_1=-\arcsin\frac{\sin^2\eta_0-\frac{1}{2}(1+\tilde\psi_2^2-\tilde\psi_1^2)}{\sqrt{\frac{1}{4}(1+\tilde\psi_2^2-\tilde\psi_1^2)^2-\tilde\psi_2^2}},
\]
and
\[
\cot^2\eta(s)=-1+\frac{1}{\frac{1}{2}\left(1+\tilde\psi_2^2-\tilde\psi_1^2\right)+D_0\sin( \pm4s\sqrt{A}+D_1)}.
\]
Integrating two first equations of the Hamiltonian system \eqref{Hsystem} gives
\[
\zeta_1(s)-\zeta_1^0=s\sqrt{A}\tilde\psi_2+\sqrt{A}\tilde\psi_1\int_0^s\tan^2\eta(s)ds
\]
\[
=s\sqrt{A}(\tilde\psi_2-\tilde\psi_1)+\int_0^s \frac{\sqrt{A}\tilde\psi_1ds}{1-\frac{1}{2}\left(1+\tilde\psi_2^2-\tilde\psi_1^2\right)- D_0\sin( \pm4s\sqrt{A}+D_1)}
\]
\[
=s\sqrt{A}(\tilde\psi_2-\tilde\psi_1)+2\arctan\frac{1}{\tilde\psi_1}\left\{\Bigg[ \frac{1}{2}\left(1-\tilde\psi_2^2+\tilde\psi_1^2\right)\Bigg] \tan\left( \pm 2s\sqrt{A}+\frac{D_1}{2}\right)- D_0\right\}
\]
\[
-2\arctan\frac{1}{\tilde\psi_1}\left\{\Bigg[ \frac{1}{2}\left(1-\tilde\psi_2^2+\tilde\psi_1^2\right)\Bigg] \tan\left(\frac{D_1}{2}\right)- D_0\right\}.
\]

Analogously,
\[
\zeta_2(s)-\zeta_2^0=s\sqrt{A}\tilde\psi_1+\sqrt{A}\tilde\psi_2\int_0^s\cot^2\eta(s)ds
\]
\[
=s\sqrt{A}(\tilde\psi_1-\tilde\psi_2)+2\arctan\frac{1}{\tilde\psi_2}\left\{\frac{1}{2}\left(1+\tilde\psi_2^2-\tilde\psi_1^2\right) \tan\left(\pm 2s\sqrt{A}+\frac{D_1}{2}\right)+ D_0\right\}
\]
\[
-2\arctan\frac{1}{\tilde\psi_2}\left\{\frac{1}{2}\left(1+\tilde\psi_2^2-\tilde\psi_1^2\right) \tan\left(\frac{D_1}{2}\right)+ D_0\right\}.
\]

Let us study solutions starting at the point $x_0=(1/\sqrt{2},0,1/\sqrt{2},0)$, or $\zeta_1^0=\zeta_2^0=0$, $\eta_0=\pi/ 4$. 
Observe that the option $\tilde\psi_1=\tilde\psi_2=0$ leads to the trivial solution $\zeta_1(s)\equiv 0$, $\zeta_2(s)\equiv 0$, and $\eta(s)=\frac{\eta-\pi/4}{t}s+\pi/4$. If $\tilde\psi_1=0$ and $\tilde\psi_2\neq 0$, then $\zeta_1(s)=s\tilde\psi_2\sqrt{A}$ and $\zeta_2(s)$, $\eta(s)$ are calculated by the above formulas.

As it has been shown in cartesian coordinates, there are infinite number of geodesics joining the point $x_0$ with a point of the vertical line.
In hyperspherical coordinates this corresponds to a point of the vertical line expressed as $\eta=\pi/4$, $\zeta_1=-\zeta_2$. Then,  $A^{(n)}=\left(\frac{\pi n}{2 t}\right)^2$, and the solution to the boundary value problem
for the functions $\zeta_1(s)$ and $\zeta_2(s)$ leads to $\tilde\psi_2-\tilde\psi_1=\frac{2\zeta_1}{\pi n}=\frac{-2\zeta_2}{\pi n}$ for $n\neq 0$. The case $n=0$ corresponds to the degenerate curve at initial point. Thus we have a countable number of geometrically different geodesics. The concrete choice of $\tilde\psi_1$ and $\tilde\psi_2$ fixes the rotation about the vertical line. 

In order to solve the boundary value problem in non-trivial situations we express implicitly $\tilde\psi_1$ and $\tilde\psi_2$ as solutions
$\tilde\psi_1=\tilde\psi_1^{(n)}(\zeta_1^0,\zeta_2^0,\eta_0, \zeta_1,\zeta_2,\eta)=\tilde\psi_1^{(n)}(0,0,\pi/4, \zeta_1,\zeta_2,\eta)$ and $\tilde\psi_2=\tilde\psi_2^{(n)}(\zeta_1^0,\zeta_2^0,\eta_0, \zeta_1,\zeta_2,\eta)=\tilde\psi_2^{(n)}(0,0,\pi/4, \zeta_1,\zeta_2,\eta)$ to the equations
\[
\zeta_1=t\sqrt{A}(\tilde\psi_2-\tilde\psi_1)+2\arctan\frac{1}{\tilde\psi_1}\left\{ \frac{1}{2}\left(1-\tilde\psi_2^2+\tilde\psi_1^2\right) \tan\left( \pm 2t\sqrt{A}+\frac{D_1}{2}\right)- D_0\right\}
\]
\[
-2\arctan\frac{1}{\tilde\psi_1}\left\{ \frac{1}{2}\left(1-\tilde\psi_2^2+\tilde\psi_1^2\right)\tan\left(\frac{D_1}{2}\right)- D_0\right\},
\]
and
\[
\zeta_2=t\sqrt{A}(\psi_1-\psi_2)+2\arctan\frac{1}{\tilde\psi_2}\left\{\frac{1}{2}\left(1+\tilde\psi_2^2-\tilde\psi_1^2\right) \tan\left( \pm 2t\sqrt{A}+\frac{D_1}{2}\right)+ D_0\right\}
\]
\[
-2\arctan\frac{1}{\tilde\psi_2}\left\{\frac{1}{2}\left(1+\tilde\psi_2^2-\tilde\psi_1^2\right) \tan\left(\frac{D_1}{2}\right)+ D_0\right\},
\]
where $A$ is a solution $A=A^{(n)}(\tilde\psi_1,\tilde\psi_2,\eta_0,\eta,t)$ given by \eqref{eqA}.

The Hamiltonian can be written as
\[
H=\frac{1}{2}A(1-(\tilde\psi_1-\tilde\psi_2)^2)=\frac{1}{2}A^{(n)}(1-(\tilde\psi_1-\tilde\psi_2)^2).
\]
If the point is not at the vertical line, then the Hamiltonian is bounded as it was shown in the cartesian coordinates. Moreover, $(1-(\tilde\psi_1-\tilde\psi_2)^2)>0$. Then, there can be only finite number of $n$  that satisfy all boundary conditions. From this formulation it is easy to see
which $n$ and which sign in \eqref{eqA}  we have to choose in order to define the minimizer. This is given by the condition of minimal value of $A^{(n)}$. In particular,
for the vertical line, $n=\pm 1$.

\subsection{Modified action}
We are aimed now at construction of the function $f$ mentioned in Section 2.
It turns out  that $\psi_1$ and $\psi_2$ are the first integrals of the Hamiltonian system \eqref{Hsystem}, which we assume to be given constants. 
The Hamiltonian system in the hyperspherical coordinates suggests the form of the modified action whereas in the cartesian coordinates
its possible form remains rather unclear. 
Let us solve \eqref{Hsystem} for the following mixed boundary conditions: 
\begin{equation}\label{bconditions}
\begin{array}{l}
\eta(0)=\eta_0,\,\, \eta(t)=\eta, \,\, \zeta_1(t)=\zeta_1, \,\,
\zeta_2(t)=\zeta_2, \bigskip \\ 
 \psi_1(0)=\psi_1, \,\, \psi_2(0)=\psi_2.
\end{array}
\end{equation}

In the classical case the action is defined on arbitrary smooth curves which join two given points. In our case the non-classical (or modified) action is defined on solutions to the above Hamiltonian system where only
the coordinate $\eta$ is given in two endpoints of the interval $[0,t]$. The coordinates $\zeta_1$ and $\zeta_2$ are given only at the right-hand endpoint~$t$.
We are looking for the modified action in the form
\[
g(\zeta_1,\zeta_2,\eta_0,\eta,\psi_1,\psi_2,t)=\psi_1(0)\zeta_1(0)+\psi_2(0)\zeta_2(0)+\int_0^t(\psi_1\dot{\zeta}_1(s)+\psi_2\dot{\zeta}_2(s)+\theta \dot{\eta}(s)-H)ds
\]
\[
=\psi_1\zeta_1+\psi_2\zeta_2+\int_0^t(\theta \dot{\eta}(s)-H)ds,
\]
inspired by the work~\cite{BealsGaveauGreiner3}. Observe also that what is written is nothing but the integration-by-parts formula.
The function $H$ does not depend on $s$ on the solutions to the Hamiltonian system (independently on boundary conditions), therefore we have
\[
g=\psi_1\zeta_1+\psi_2\zeta_2-\frac{t}{2}((\psi_1\tan\eta_0+\psi_2\cotan\eta_0)^2+4\theta_0^2)
+\int_0^t4\theta^2(s)ds,
\]
where $\theta_0$ is the constant of integration.
In order to calculate our modified action we have to 
\begin{itemize}
\item calculate the integral;
\item represent $\theta_0$, or equivalently $A$ in terms of $\eta_0$, $\eta$, $\psi_1$, $\psi_2$.
\end{itemize} 

Let us remark here that we do not specify at the moment the value of $\theta_0$ (or $A$) from the countable number of possible values given by \eqref{eqA}. This will be done later in Section~\ref{Sec8}.

Now let us calculate the integral in the modified action.
The action $g$ admits the form
\[
g=\psi_1\zeta_1+\psi_2\zeta_2+\frac{t}{2}((\psi_1\tan\eta_0+\psi_2\cot\eta_0)^2+4\theta_0^2)-2t\psi_1\psi_2
-\int_0^t (\psi_1^2\tan^2\eta(s)+\psi_2^2\cot^2\eta(s))ds
\]
\[
=\psi_1\zeta_1+\psi_2\zeta_2+\frac{t}{2}((\psi_1\tan\eta_0+\psi_2\cot\eta_0)^2+4\theta_0^2)+t(\psi_1-\psi_2)^2
\]
\[
-\int_0^t \left(\frac{\psi_1^2}{1-\frac{1}{2}\left(1+\frac{\psi_2^2-\psi_1^2}{A}\right) -D_0\sin(4s\sqrt{A}+D_1)}+\frac{\psi_2^2}{\frac{1}{2}\left(1+\frac{\psi_2^2-\psi_1^2}{A}\right) +D_0\sin(4s\sqrt{A}+D_1)}\right)ds.
\]
We observe that
\[
\frac{1}{4}\left(1+\frac{\psi_2^2-\psi_1^2}{A}\right)^2-D_0^2=\frac{\psi_2^2}{A}>0,
\]
and
\[
\left(1-\frac{1}{2}\left(1+\frac{\psi_2^2-\psi_1^2}{A}\right)\right)^2-D_0^2=\frac{\psi_1^2}{A}>0.
\]
Thus, integrating gives us
\[
g=\psi_1\zeta_1+\psi_2\zeta_2+\frac{t}{2}A+\frac{t}{2}(\psi_1-\psi_2)^2
\]
\[
-2\psi_1\arctan\frac{\sqrt{A}}{\psi_1}\left\{\Bigg[ 1-\frac{1}{2}\left(1+\frac{\psi_2^2-\psi_1^2}{A}\right)\Bigg] \tan\left(2t\sqrt{A}+\frac{D_1}{2}\right) -D_0\right\}
\]
\[
+2\psi_1\arctan\frac{\sqrt{A}}{\psi_1}\left\{\Bigg[ 1-\frac{1}{2}\left(1+\frac{\psi_2^2-\psi_1^2}{A}\right)\Bigg] \tan\frac{D_1}{2} -D_0\right\}
\]
\[
-2\psi_2\arctan\frac{\sqrt{A}}{\psi_2}\left\{\frac{1}{2}\left(1+\frac{\psi_2^2-\psi_1^2}{A}\right) \tan\left(2t\sqrt{A}+\frac{D_1}{2}\right) +D_0\right\}
\]
\[
+2\psi_2\arctan\frac{\sqrt{A}}{\psi_2}\left\{\frac{1}{2}\left(1+\frac{\psi_2^2-\psi_1^2}{A}\right) \tan\frac{D_1}{2} +D_0\right\}.
\]
The important particular case which we shall use in forthcoming sections is $\psi_2=-\psi_1$, that leads to a much simpler formula
\[
g(\zeta_1,\zeta_2,\pi/4,\eta,\psi_1,-\psi_1,t)=\psi_1(\zeta_1-\zeta_2)+\frac{At}{2}+2t\psi_1^2-2\psi_1\arctan\frac{\psi_1}{\sqrt{A}}\tan 4t\sqrt{A}.
\]

\section{Generalized Hamilton-Jacobi equation}\label{Sec7}

\begin{theorem} Let $\zeta_1$, $\zeta_2$, and $\eta$ be fixed, and let $\zeta_1(s)=\zeta_1(s;\psi_1,\psi_2,\zeta_1,\zeta_2, \eta,\eta_0, t)$ $\zeta_2(s)=\zeta_2(s;\psi_1,\psi_2,\zeta_1,\zeta_2, \eta,\eta_0, t)$, and $\eta(s)=\eta(s;\psi_1,\psi_2, \eta,\eta_0, t)$ be solutions to the Hamiltonian system \eqref{Hsystem} with the mixed boundary conditions \eqref{bconditions} and with the Hamiltonian $H(\zeta_1(s),\zeta_2(s),\eta(s),\psi_1,\psi_2,\theta(s))$. The modified action
\[
g=\psi_1\zeta_1+\psi_2\zeta_2+\int_0^t(\theta \dot{\eta}(s)-H)ds
\]
satisfies the Hamilton-Jacobi equation
\[
\frac{\partial g}{\partial t}+H(\zeta_1,\zeta_2,\eta,\nabla g)=0,
\]
where the gradient $\nabla$ is taken with respect to the coordinates of the endpoint $(\zeta_1,\zeta_2,\eta)$.
\end{theorem}
\begin{proof}
Let us calculate the derivatives of $g$ with respect to $\zeta_1$, $\zeta_2$, and $\eta$ explicitly.
\[
\frac{\partial g}{\partial \zeta_1}=\psi_1+\int_0^t\left(\frac{\partial \theta(s)}{\partial \zeta_1}\dot{\eta}+\theta(s) \frac{d}{ds}\frac{\partial \eta(s)}{\partial \zeta_1}-\frac{\partial H}{\partial \eta(s)}\frac{\partial \eta(s)}{\partial \zeta_1}- \frac{\partial H}{\partial \theta(s)}\frac{\partial \theta(s)}{\partial \zeta_1}\right)ds
\]
\[
=\psi_1+\int_0^t\left(\frac{d}{ds}\Big[\theta(s)\frac{\partial \eta(s)}{\partial \zeta_1}\Big]\right)ds=\psi_1.
\]
We used that $\frac{\partial H}{\partial \theta(s)}=\dot{\eta}$,  $\frac{\partial H}{\partial \eta(s)}=-\dot{\theta}$. Moreover, $\eta(s)$ does not depend on $\zeta_1$ for $s\in[0,t]$. Here $\dot{\eta}$ means the derivative with respect to $s$. Analogously,
\[
\frac{\partial g}{\partial \zeta_2}=\psi_2.
\]
We continue with
\[
\frac{\partial g}{\partial \eta}=\int_0^t\left(\frac{\partial \theta(s)}{\partial \eta}\dot{\eta}+\theta(s) \frac{d}{ds}\frac{\partial \eta(s)}{\partial \eta}-\frac{\partial H}{\partial \eta(s)}\frac{\partial \eta(s)}{\partial \eta}- \frac{\partial H}{\partial \theta(s)}\frac{\partial \theta(s)}{\partial \eta}\right)ds
\]
\[
=\int_0^t\left(\frac{d}{ds}\Big[\theta(s)\frac{\partial \eta(s)}{\partial \eta}\Big]\right)ds=\theta(s)\frac{\partial \eta(s)}{\partial \eta}\Bigg|_{s=0}^{s=t}.
\]
The latter derivative is taken with respect to $t$ as
\[
\frac{\partial g}{\partial t}=(\theta(s)\dot{\eta}-H)_{s=t}+\int_0^t\left(\frac{\partial \theta(s)}{\partial t}\dot{\eta}+\theta(s)\frac{\partial }{\partial t}\dot{\eta}(s)-\frac{\partial H}{\partial \eta(s)}\frac{\partial \eta(s)}{\partial t}-\frac{\partial H}{\partial \theta(s)}\frac{\partial \theta(s)}{\partial t}\right)ds
\]
\[
=(\theta(s)\dot{\eta}-H)_{s=t}+\theta(s)\frac{\partial \eta(s)}{\partial t}\Bigg|_{s=0}^{s=t}.
\]
We need to calculate
\[
\frac{\partial \eta(s)}{\partial \eta}\Bigg|_{s=0}^{s=t} \quad\mbox{and \ }\frac{\partial \eta(s)}{\partial t}\Bigg|_{s=0}^{s=t}.
\]
Let us rewrite the equation \eqref{eqeta} in the form
\[
P(\eta(s),A)=4s\sqrt{A},
\]
where $A=A^{(n)}(\psi_1,\psi_2,\eta_0,\eta,t)$ satisfies the equation
\[
P(\eta(t),A)=4t\sqrt{A}.
\]
Then
\[
\frac{d P(\eta(s), A)}{d \eta}  =  \frac{\partial P(\eta(s), A)}{\partial\eta(s)}\frac{\partial \eta(s)}{\partial \eta}+\frac{\partial P(\eta(s), A)}{\partial A}\frac{\partial A}{\partial \eta}=\frac{2s}{\sqrt{A}}\frac{\partial A}{\partial \eta}.
\]
Substituting first $s=t$, and then differentiating with respect to $\eta$ we come to
\[
\frac{d P(\eta, A)}{d\eta} =  \frac{\partial P(\eta, A)}{\partial\eta}+\frac{\partial P(\eta, A)}{\partial A}\frac{\partial A}{\partial \eta}=\frac{2t}{\sqrt{A}}\frac{\partial A}{\partial \eta}.
\]
Substituting $s=t$ in the first equation and subtracting the second one we obtain
\[
\frac{\partial \eta(s)}{\partial \eta}\Big|_{s=t}=1.
\]
Substituting $s=0$ in the first equation we get
\[
\frac{\partial \eta(s)}{\partial \eta}\Big|_{s=0}=-\frac{\partial P(\eta(s), A)}{\partial A}\frac{\partial A}{\partial \eta}\left( \frac{\partial P(\eta(s), A)}{\partial\eta(s)}\right)^{-1}\Bigg|_{s=0}=0.
\]
Analogously, calculating the derivative
\[
\frac{\partial P(\eta(s), A)}{\partial t},
\]
and substituting $s=0$ we get $\frac{\partial \eta(s)}{\partial t}\Big|_{s=0}=0$. Finally, we have
\[
0=\frac{\eta(t; \psi_1,\psi_2,\eta_0,\eta, t)}{\partial t}=\left(\dot{\eta}+\frac{\eta(s; \psi_1,\psi_2,\eta_0,\eta, t)}{\partial t}\right)_{s=t},
\]
which implies  the relation $\frac{\partial g}{\partial t}=-H(\zeta_1,\zeta_2,\eta,\psi_1,\psi_2,\theta)$.
Now we are able to calculate $H(\zeta_1,\zeta_2,\eta,\nabla g)$ as
\[
H(\zeta_1,\zeta_2,\eta,\nabla g)=\frac{1}{2}\Big[\left(\frac{\partial g}{\partial \zeta_1}\right)^2\tan^2\eta+\left(\frac{\partial g}{\partial \zeta_2}\right)^2\cot^2\eta+2\frac{\partial g}{\partial \zeta_1}\frac{\partial g}{\partial \zeta_2}+4\left(\frac{\partial g}{\partial \eta}\right)^2\Big]
\]
\[
=\frac{1}{2}(\psi_1^2\tan^2\eta+\psi_2^2\cot^2\eta+2\psi_1\psi_2+4\theta^2)_{s=t}=H(\zeta_1,\zeta_2,\eta,\psi_1,\psi_2,\theta)=-\frac{\partial g}{\partial t}.
\]
This finishes the proof.
\end{proof}

Let us introduce the differential operator $T_g=\frac{\partial
}{\partial t}+\nabla_0g\cdot\nabla_0$, where $g$ is the modified
action and $\nabla_0g$ is the horizontal gradient $(Xg,Yg)$. The operator $T$ was
deeply studied by Beals, Gaveau and Greiner, see e.g.,
in~\cite{BealsGaveauGreiner2}. Let
$\gamma(s)=(\zeta_1(s),\zeta_2(s),\eta(s))$, $s\in[0,t]$, be a
geodesic, that is a projection of a solution of the Hamiltonian
system onto the underlying manifold. Since $\frac{\partial
g}{\partial \zeta_i}=\psi_1$, $i=1,2$, $\frac{\partial g}{\partial
\eta}=\theta$, by definition, we have  that the formal symbol of the vector fields $X,Y$
defines the horizontal gradient of $g$ along the geodesic:
\begin{eqnarray*}\big(X(\gamma(s),\psi,
\theta(s)),Y(\gamma(s),\psi_i, \theta(s))\big) & = &
\big(X(\gamma(s),\frac{\partial g}{\partial
\zeta_i},\frac{\partial g}{\partial
\eta}),Y(\gamma(s),\frac{\partial g}{\partial
\zeta_i},\frac{\partial g}{\partial \eta})\big) \\ & = &
\nabla_0g(\gamma(s)).\end{eqnarray*} The projection onto the
manifold defines the horizontal gradient of $\gamma$:
$$\big(X(\gamma),Y(\gamma)\big)=\nabla_0(\gamma(s)).$$

For any function $h(\zeta_i,\eta,t)$  let us define its value
 along a geodesic
$\gamma(s)=(\zeta_i(s),\eta(s))$ by
$h(s)=h(\zeta_i(s),\eta(s),s)$, $s\in[0,t]$. Then the differential
operator $T_g(h(\gamma))$ is the differentiation along the geodesic
$\gamma$. Indeed, $$T_g(h(\gamma))=\frac{\partial h}{\partial
t}(s)+\nabla_0h(\zeta_i,\eta)\cdot\nabla_0(\gamma)(s)
=\frac{\partial h}{\partial
t}(s)+\nabla_0h(\zeta_i,\eta)\cdot\nabla_0(g(\gamma))(s).$$

Denote $\nabla
h=\big(\frac{\partial h}{\partial \zeta_1},\frac{\partial h}{\partial \zeta_2},\frac{\partial
h}{\partial \eta}\big)$,.

\begin{proposition}\label{prop01} If a function
$h(\zeta_i,\eta,t)$ satisfies the Hamilton--Jacobi equation
$h_t=-H(\nabla h)$, then its derivative
$h_t(\zeta_i,\eta,t)=\frac{\partial h(\zeta_i,\eta,t)}{\partial
t}$ possesses the properties:
$$T_h(h_t)=\frac{\partial h_t}{\partial
t}+\nabla_0h\cdot\nabla_0h_t=0,$$ and 
\begin{equation}\label{T(h)}T_h(h)=-h_t.\end{equation}
\end{proposition}

\begin{proof}
If the Hamilton - Jacobi equation is satisfied for $h$, then
$$T_h(h_t)=\frac{\partial h_t}{\partial t}+\nabla_0 h\cdot\nabla_0
h_t=\frac{\partial h_t}{\partial t}+\frac{1}{2}\frac{\partial
}{\partial t}\Big(\nabla_0 h\cdot\nabla_0 h\Big)=\frac{\partial
}{\partial t}\Big(h_t+H(\nabla h)\Big)=0.$$ Applying the operator
$T$ to $h$ we obtain
$$T_h(h)=\frac{\partial h}{\partial t}+\nabla_0 h\cdot\nabla_0 h
= -H(\nabla h)+2H(\nabla h)=H(\nabla h)=-h_t.$$
\end{proof}

Observe that the modified action satisfies the stretching property
\[
g(\zeta_1,\zeta_2,\eta_0,\eta,\psi_1,\psi_2,t)=\lambda g(\zeta_1,\zeta_2,\eta_0,\eta,\frac{\psi_1}{\lambda},\frac{\psi_2}{\lambda},\lambda t),
\]
and construct the function
\[
f(\zeta_1,\zeta_2,\eta_0,\eta,\psi_1,\psi_2)=g(\zeta_1,\zeta_2,\eta_0,\eta,\psi_1,\psi_2,1).
\]
\begin{theorem}\label{ThHJ}
The function $f$ is a solution to the generalized Hamilton--Jacobi equation
\[
\psi_1\frac{\partial f}{\partial \psi_1}+\psi_2\frac{\partial f}{\partial \psi_2}+H(\zeta_1,\zeta_2,\eta,\nabla f)=f.
\]
\end{theorem}
\begin{proof}
By the above stretching property of the function $g$ we have
\[
g(\zeta_1,\zeta_2,\eta_0,\eta,\psi_1,\psi_2,t)=\frac{1}{t}g(\zeta_1,\zeta_2,\eta_0,\eta,t\psi_1,t\psi_2,1)=\frac{1}{t}f(\zeta_1,\zeta_2,\eta_0,\eta,t\psi_1,t\psi_2).
\]
Therefore,
\[
\frac{\partial g}{\partial t}=-\frac{1}{t^2}f+\frac{1}{t}\left(\psi_1\frac{\partial f}{\partial \psi_1}+\psi_2\frac{\partial f}{\partial \psi_2}\right),
\]
for any fixed $t$, in particular for $t=1$.
On the other hand $g$ satisfies the Hamilton--Jacobi equation. This finishes the proof.
\end{proof}

\section{Modified action as a distance function}\label{Sec8}

Theorem \ref{ThHJ} implies that if $\psi_1$ and $\psi_2$ are critical points for the modified action $f$, then $f$ is equal to the Hamiltonian $H(\zeta_1,\zeta_2,\eta,\nabla f)$. This allows us to interpret $f$ as a distance function.

\begin{theorem} Suppose that the endpoint $\zeta_1$, $\zeta_2$ and $\eta$ does not belong to the vertical line passing through the initial point $P=(0,0,\eta_0)$.
Among all critical points $\psi_1=\psi_1^{(n)}(\zeta_1,\zeta_2,\eta_0,\eta)$ and $\psi_2=\psi_2^{(n)}(\zeta_1,\zeta_2,\eta_0,\eta)$ of the  modified action $f(\zeta_1,\zeta_2,\eta_0,\eta,\psi_1,\psi_2)$, i.e., those points
satisfying the equations $\partial f/\partial\psi_1=0$ and $\partial f/\partial\psi_2=0$, there exist exactly one such that the action $f$  evaluated at this critical point is a distance function from the point $P=(0,0,\eta_0)$ to $Q=(\zeta_1,\zeta_2,\eta)$. Every such critical value
$\psi_1^{(n)}$ and $\psi_2^{(n)}$ defines a unique geodesic joining $P$ and $Q$.
\end{theorem}
\begin{proof}
Let us calculate the derivative
\[
\frac{\partial f}{\partial \psi_1}=\zeta_1+\int_0^1\left(\frac{\partial \theta(s)}{\partial \psi_1}\dot{\eta}+\theta(s) \frac{d}{ds}\frac{\partial \eta(s)}{\partial \psi_1}-\frac{\partial H}{\psi_1}-\frac{\partial H}{\partial \eta(s)}\frac{\partial \eta(s)}{\partial \psi_1}- \frac{\partial H}{\partial \theta(s)}\frac{\partial \theta(s)}{\partial \psi_1}\right)ds
\]
\[
=\zeta_1+\theta(s)\frac{\partial \eta(s)}{\partial \psi_1}\Bigg|_{s=0}^{s=t}-\int_0^t\frac{\partial H}{\partial \psi_1}ds.
\]
We need to calculate
\[
\frac{\partial \eta(s)}{\partial \psi_1}\Bigg|_{s=0}\quad\mbox{and \ }\frac{\partial \eta(s)}{\partial \psi_1}\Bigg|_{s=t}.
\]
Similarly to the previous section we use the derivatives of the function $P(\eta(s),A)$ as
\[
\frac{\partial P(\eta(s), A)}{\partial\psi_1}  =  \frac{\partial P(\eta(s), A)}{\partial \psi_1}+  \frac{\partial P(\eta(s), A)}{\partial\eta(s)}\frac{\partial \eta(s)}{\partial \psi_1}+\frac{\partial P(\eta(s), A)}{\partial A}\frac{\partial A}{\partial \psi_1}=\frac{2s}{\sqrt{A}}\frac{\partial A}{\partial \psi_1}.
\]
Substituting $s=0$ in the latter equation we get
\[
\frac{\partial \eta(s)}{\partial \psi_1}\Bigg|_{s=0}=0.
\]
Consider now the derivative
\[
\frac{\partial P(\eta, A)}{\partial\psi_1}  =  \frac{\partial P(\eta, A)}{\partial \psi_1}+\frac{\partial P(\eta, A)}{\partial A}\frac{\partial A}{\partial \psi_1}=\frac{2t}{\sqrt{A}}\frac{\partial A}{\partial \psi_1}.
\]
Then
\[
\frac{\partial P(\eta(s), A)}{\partial\psi_1}\Bigg|_{s=t} -\frac{\partial P(\eta, A)}{\partial\psi_1} = \frac{\partial P(\eta(s), A)}{\partial\eta(s)}\frac{\partial \eta(s)}{\partial \psi_1}\Bigg|_{s=t}=0.
\]
Hense,
\[
\frac{\partial \eta(s)}{\partial \psi_1}\Bigg|_{s=t}=0.
\]
Therefore,
\begin{equation}\label{crit1}
\frac{\partial f}{\partial \psi_1}=\zeta_1-\int_0^1\dot{\zeta}_1ds=\zeta_1^0=0.
\end{equation}
Analogously,
\begin{equation}\label{crit2}
\frac{\partial f}{\partial \psi_2}=\zeta_2-\int_0^1\dot{\zeta}_2ds=\zeta_2^0=0.
\end{equation}
On the other hand, the equations \eqref{crit1} and \eqref{crit2} can be rewritten in terms of the function $\eta(s)$ as
\[
\frac{\partial f}{\partial \psi_1}=\zeta_1-\int_0^1(\psi_1 \tan^2\eta(s)+\psi_2)ds=0,\quad \frac{\partial f}{\partial \psi_2}=\zeta_2-\int_0^1(\psi_2 \cot^2\eta(s)+\psi_1)ds=0.
\]
The critical points $\psi_1=\psi_1^{(n)}(\zeta_1,\zeta_2,\eta_0,\eta)$ and $\psi_2=\psi_2^{(n)}(\zeta_1,\zeta_2,\eta_0,\eta)$ are defined as solutions to the above system of two equations, which is the
same as the system which defines reparametrization of geodesics  for $\zeta_1^0=\zeta_2^0=0$. Thus, every critical point $\psi_1=\psi_1^{(n)}(\zeta_1,\zeta_2,\eta_0,\eta)$ and $\psi_2=\psi_2^{(n)}(\zeta_1,\zeta_2,\eta_0,\eta)$ defines a unique geodesic starting at the point $P$ and ending at the point $Q$ as stated in the theorem. Since the point $Q$ does not belong to the vertical line, there is a finite number
of geodesics joining $P$ and $Q$, see Theorem \ref{numbergeod2}. We choose the geodesic which minimizes the Hamiltonian ($C^2+D^2$ in terms of Theorem \ref{numbergeod2}). Thus, the function $f$
evaluated in its critical points corresponding to this geodesic satisfies all the properties of a distance, i.e., non-negative, vanishes
if and only if $P=Q$, and satisfies the triangle inequality.
\end{proof}

We consider now the modified action $f$ restricted to the characteristic variety in the cotangent bundle
defined by $\psi_1=\tau$, $\psi_2=-\tau\tan^2\eta_0$, $\theta=0$. This variety is a singular set because
the Hamiltonian $H$ vanishes there. The modified action does not depend on the variable $\theta$ of the tangent bundle. Choosing the initial point $(0,0,\pi/4)$ in the phase manifold with coordinates
$(\zeta_1,\zeta_2,\eta)$, we set $\psi_1=\tau$, $\psi_2=-\tau$.  Let us suppose the sign (+) in all formulas for the modified action. This simplify significantly the analytic expression for $f$ which becomes of the form
\[
f(\psi_1,\psi_2,\zeta_1,\zeta_2,\eta, \eta_0)=f(\tau,-\tau,\zeta_1,\zeta_2,\eta, \pi/4)
\]
\[
=\tau(\zeta_1-\zeta_2)+\frac{1}{2}A+2\tau^2-2\tau \arctan\left(4\frac{\tau}{\sqrt{A}}\frac{\tan 2\sqrt{A}}{1-\tan^22\sqrt{A}}\right)
\]
\[
=\tau(\zeta_1-\zeta_2)+\frac{1}{2}A+2\tau^2-2\tau \arctan\left(\frac{2\tau}{\sqrt{A}}\tan 4\sqrt{A}\right),
\]
where $A>4\tau^2$, and satisfies the equation
\[
\sin^2\eta=\frac{1}{2}+(\sin 4\sqrt{A})\sqrt{\frac{1}{4}-\frac{\tau^2}{A}}.
\]
Among the infinite number of the solutions to this equation we choose the first one. It gives the minimum to the Hamiltonian.
We see that the action $f(\tau,-\tau,\zeta_1,\zeta_2,\eta, \pi/4)$ has singularities at all points
\[
\tau=\pm\frac{1}{16}(\pi+2\pi n)\sin 2\eta.
\]
However, the term $\tau^2$ is dominating and the function $\exp(-f)$ exponentially decays
as $\tau\to\pm \infty$ and integrable. This contrasts the case of non-compact sub-Riemannian
manifolds considered earlier in, e.g.,  \cite{BealsGaveauGreiner1, BealsGaveauGreiner2, BealsGaveauGreiner3, CalinChangMarkina2}, where the authors had to avoid non-integrable singularities introducing complex modified action.

\section{Volume element}

The aim of this section is to deduce the transport equation for the volume element that serves for the heat kernel, and to find  the fundumental solution to the sub-Laplacian equation in the case of 3-D sub-Riemannian sphere.

\subsection {Volume element for the heat kernel.} 

 We shall use the notation of  Sections \ref{Sec3} and~\ref{Modif.function}. Choose the point $(0,0,\frac{\pi}{4})$ as the center for the heat kernel just to have  simpler formulas. The characteristic variety is the straight line given by
$$\left\{\begin{array}{lll}
\psi_1 & = & \tau\cotan\eta\\
\psi_2 & = & -\tau\tan\eta\\
\theta & = & 0
\end{array}
\right.$$ 
We start from the volume element $v$ for the heat kernel. Let $f$ be the modified action  studied at Section~
\ref{Modif.function}, ~\ref{Sec7}, ~\ref{Sec8}. We remind that it can be considered as the
square of the distance from the point $(0,0,\eta_0)$, and in
particular, from $\eta_0=\frac{\pi}{4}$ to some point
$(\zeta_1,\zeta_2,\eta)$. Let us integrate the function
$e^{-\frac{f}{u}}$ over the characteristic variety at
$(0,0,\frac{\pi}{4})$ with respect to a measure
$v(\zeta_1,\zeta_2,\eta, \tau)$. The characteristic variety at
$(0,0,\frac{\pi}{4})$ gives us the relation between $\psi_1$ and
$\psi_2$, namely, $\psi_1=-\psi_2=\tau$. Since $f$ satisfies the
generalized Hamilton--Jacobi equation for any $\psi_1,\psi_2$,  we
have the equation \begin{equation}\label{genHJ}\tau\frac{d
f}{d\tau}+H(\zeta_1,\zeta_2,\eta,\tau,\nabla
f)=f(\zeta_1,\zeta_2,\eta,\tau),\end{equation} where we write
$\tau\frac{df}{d\tau}=\Big(\psi_1\frac{\partial
f}{\partial\psi_1}+\psi_2\frac{\partial
f}{\partial\psi_2}\Big)\vert_{\psi_1=-\psi_2=\tau}$. Let us look at
some interesting properties of  $f(\tau)$. Define
$f(\zeta_i,\eta,\tau,-\tau)=\tau h(\zeta_i,\eta,\tau)$, $i=1,2$
and denote $f_{\tau}=\frac{df(\tau,-\tau)}{d\tau}$, $h_\tau=\frac{\partial
h}{\partial \tau}$, $\nabla f=(\frac{\partial f}{\partial
\zeta_1},\frac{\partial f}{\partial \zeta_2},\frac{\partial
f}{\partial \eta})$, $\nabla h=(\frac{\partial h}{\partial
\zeta_1},\frac{\partial h}{\partial \zeta_2},\frac{\partial
h}{\partial \eta})$.

\begin{proposition}\label{pr}
Under the above mentioned notations the following properties hold:
\begin{itemize}
\item[(i)]{the generalized Hamilton-Jacobi equation~\eqref{genHJ} for $f$ implies that $h$ satisfies the Hamilton-Jacobi equation;}
\item[(ii)]{$f_{\tau}=-\tau H(\nabla h)+h$;}
\item[(iii)]{$f_{\tau}=\tau h_{\tau}+h=0$;}
\item[(iv)]{$T_h(f_{\tau})=\frac{d\, f_{\tau}}{d\tau}+\nabla_0h\cdot \nabla_0f_{\tau}=0$;}
\item[(v)]{$T_h(f)=\frac{d f}{d\tau}+\nabla_0h\cdot \nabla_0f=h-\tau h$.}
\end{itemize}
\end{proposition}

\begin{proof}
In order to prove (i) we have to show that $\frac{\partial
h}{\partial\tau}+H(\nabla h)=0$. Indeed,
$$\tau\frac{df}{d\tau}-f=-H(\nabla f)\quad\Rightarrow\quad
\tau(h+\tau\frac{\partial h}{\partial\tau})-\tau h
=-\tau^2H(\nabla h)\quad\Rightarrow\quad\frac{\partial h}{\partial\tau}+H(\nabla h)=0,$$
since the Hamiltonian is a homogeneous function of  order $2$ with respect to $\nabla h$.
In order to prove  (ii), we write
$$f_{\tau}=\frac{1}{\tau}\big(-H(\nabla f)+f\big)=-\tau H(\nabla h)+
h.$$ Then,  $H(\nabla h)=-h_{\tau}$ implies (iii).

Observe that  the derivative $\frac{\partial f}{\partial\tau}$ of the operator $T_h$ in (iv) is understood as
a complete derivative of $f$ with respect to $\tau$:
$\frac{df(\tau,-\tau)}{d\tau}$. Moreover $\frac{d
h}{d\tau}+\nabla_0h\cdot \nabla_0h=T_hh$, and if $h$ satisfies the
Hamilton--Jacobi equation, then $T_h(\frac{\partial
h}{\partial\tau})=0$ and $T_h(h)=-h_{\tau}$ by
Proposition~\ref{prop01}. Thus,
\begin{equation*}T_h(f_{\tau})=T_h(\tau h_{\tau}+h)  =  h_{\tau}+\tau
T_h(h_{\tau})+T_h(h)=0.
\end{equation*}

In order to prove (v) we exploit the Hamilton--Jacobi equation for $h$. Just calculate
\begin{equation*}
T_h(f)=\frac{d f}{d\tau}+\nabla_0 h\cdot \nabla_0 f  =  h+\tau
h_{\tau}+2\tau H({\nabla h})=h+\tau h_{\tau}-2\tau h_{\tau}=h-\tau
h_{\tau}.
\end{equation*}
\end{proof}

For  simplicity  we follow the
ideas of~\cite{BealsGaveauGreiner11,Greiner}, and assume that the
volume element $v$ is represented as $v=V\frac{d f}{d \tau}=VE,$ where
$E=\frac{d f}{d \tau}$. Then, let us look for the heat kernel in the form
\begin{equation}\label{eq:heatkernel}P_u(\zeta_1,\zeta_2,\eta)=\frac{C}{u^q}\int_{-\infty}^{+\infty}
e^{-\frac{f(\zeta_1,\zeta_2,\eta,\tau)}{u}}V(\zeta_1,\zeta_2,\eta,\tau)E(\zeta_1,\zeta_2,\eta,\tau)\,d\tau.\end{equation}
Let us deduce the transport equation for
$v(\zeta_1,\zeta_2,\eta,\tau)$. The relation
$$(\Delta_X-\frac{\partial}{\partial
u})P_u(\zeta_1,\zeta_2,\eta)=C\int_{-\infty}^{+\infty}
(\Delta_X-\frac{\partial}{\partial
u})\Big(u^{-q}e^{-\frac{f(\zeta_1,\zeta_2,\eta,\tau)}{u}}(VE)(\zeta_1,\zeta_2,\eta,\tau)\Big)\,d\tau,$$
is to be satisfied. One calculates
\begin{eqnarray}\label{eq:1}
\frac{\partial}{\partial u}\Big(u^{-q}e^{-\frac{f(\zeta_1,\zeta_2,\eta,\tau)}{u}}(VE)(\zeta_1,\zeta_2,\eta,\tau)\Big)=
\Big(\frac{e^{-\frac{f}{u}}(VE)}{u^{q+1}}\Big)
\Big(-q+\frac{f}{u}\Big),
\end{eqnarray}
and
\begin{eqnarray}\label{eq:2}
\Delta_X\Big(u^{-q}e^{-\frac{f}{u}}(VE)\Big) & = &
\frac{e^{-\frac{f}{u}}}{u^{q+1}}\Big(-(VE) \Delta_X f-\nabla_0
f\cdot \nabla_0 (VE) \nonumber
\\ & + & \frac{(VE)}{u}\big(f-\tau\frac{d f}{d\tau}\big)+u\Delta_X (VE)\Big),
\end{eqnarray} where  $\nabla_0 f=(X f,Y f)$ stands for the horizontal gradient. The Hamilton--Jacobi equation yelds
$$\frac{1}{2}|\nabla_0 f|^2=H(\nabla f)=f-\tau\frac{d f}{d\tau}.$$
Subtracting~\eqref{eq:1} from~\eqref{eq:2} we arrive at
\begin{eqnarray*}
(\Delta_X-\frac{\partial}{\partial
u})\Big(\frac{e^{-\frac{f}{u}}(VE)}{u^{q}}\Big) & = &
\Big(\frac{e^{-\frac{f}{u}}}{u^{q+1}}\Big)\Big(\big(q-\Delta_X
f\big)(VE) - \nabla_0 f\cdot \nabla_0 (VE) \\ & - &
\frac{1}{u}(VE)\tau \frac{d f}{d\tau}
 + u\Delta_X(VE)\Big). \end{eqnarray*}
Substituting the term $-\frac{1}{u}e^{-\frac{f}{u}}(VE)\tau
\frac{d f}{d\tau}$ from the equality
$$-\frac{1}{u}e^{-\frac{f}{u}}(VE)\tau \frac{d
f}{d\tau}=\frac{\partial}{\partial
\tau}\Big(e^{-\frac{f}{u}}(VE)\tau \Big)-(VE)e^{-\frac{f}{u}}-\tau
e^{-\frac{f}{u}}\frac{\partial (VE)}{\partial \tau},$$ one deduces
\begin{eqnarray}\label{heq1}(\Delta_X-\frac{\partial}{\partial u})\Big(\frac{e^{-\frac{f}{u}}(VE)}{u^{q}}\Big) & = &
\Big(\frac{e^{-\frac{f}{u}}}{u^{q+1}}\Big)\Big(\big(q-1-\Delta_X
f\big)(VE)-\nabla_0f\cdot \nabla_0(VE)-\tau\frac{\partial
(VE)}{\partial\tau}\nonumber\\ & + & u\Delta_X(VE)\Big)
+\frac{1}{u^{q+1}}\frac{\partial}{\partial\tau}\Big(e^{-\frac{f}{u}}\tau
(VE)\Big).\end{eqnarray}

Let us write $f(\tau)=\tau h(\tau)$. Set
$T_h=\frac{\partial}{\partial \tau}+\nabla_0 h\cdot\nabla_0$. The
derivative $\frac{\partial f}{\partial \tau}$ is thought of
as the derivative $\frac{df(\tau,-\tau)}{d\tau}=f_{\tau}$. In more condensed form we write $$-\nabla_0f\cdot \nabla_0(VE)-\tau\frac{\partial
(VE)}{\partial\tau}=-\tau\Big(\frac{\partial
(VE)}{\partial\tau}+\nabla_0h\cdot \nabla_0(VE)\Big)=-\tau
T_h(VE).$$ Substituting $E=f_{\tau}$ we get $$-\tau
T_h(VE)=-\tau\Big(f_{\tau}T_h(V)+VT_h(f_{\tau})\Big)=-\tau
f_{\tau}T_h(V)$$ by (iv). Now~\eqref{heq1} takes the form
\begin{eqnarray}\label{heq2}(\Delta_X-\frac{\partial}{\partial u})\Big(\frac{e^{-\frac{f}{u}}(Vf_{\tau})}{u^{q}}
\Big) & = &
\Big(\frac{e^{-\frac{f}{u}}}{u^{q+1}}\Big)\Big(f_{\tau}\Big[\big(q-1-\Delta_X
f\big)(V)-\tau T_h(V)\Big]+ u\Delta_X(VE)\Big)\nonumber
\\ & + &\frac{1}{u^{q+1}}\frac{\partial}{\partial\tau}\Big(e^{-\frac{f}{u}}\tau
(VE)\Big).\end{eqnarray} In the following step we observe that
\begin{eqnarray*}&   \frac{e^{-\frac{f}{u}}}{u^{q+1}}\,  f_{\tau} \,\Big[\big(q-1-\Delta_X
f\big)(V)-\tau T_h(V)\Big] \\  = &
-\frac{\partial}{\partial\tau}\Big[\frac{e^{-\frac{f}{u}}}{u^{q}}\Big(\big(q-1-\Delta_X
f\big)(V)-\tau T_h(V)\Big)\Big] \\  + &
\Big(\frac{e^{-\frac{f}{u}}}{u^{q}}\Big)\frac{\partial}{\partial\tau}\Big(\big(q-1-\Delta_X
f\big)(V)-\tau T_h(V)\Big), \end{eqnarray*} and substitute the result
into~\eqref{heq2}. This leads to
\begin{eqnarray}\label{heq3}(\Delta_X-\frac{\partial}{\partial u})\Big(\frac{e^{-\frac{f}{u}}(Vf_{\tau})}{u^{q}}
\Big) & = &
\Big(\frac{e^{-\frac{f}{u}}}{u^{q}}\Big)\Big(\frac{\partial}{\partial\tau}\Big[\big(q-1-\Delta_X
f\big)(V)-\tau T_h(V)\Big]+ \Delta_X(Vf_{\tau})\Big)\nonumber
\\ & + &\frac{1}{u^{q}}\frac{\partial}{\partial\tau}\Big(\big(e^{-\frac{f}{u}}\big)
\big(\tau T_h(V)+(\Delta_X f-q+1)V\big)\Big)\\ & + &
\frac{1}{u^{q+1}}\frac{\partial}{\partial\tau}\Big(e^{-\frac{f}{u}}\tau
f_{\tau}V\Big).\nonumber\end{eqnarray}

Working with the first term one obtains
\begin{eqnarray*}
 \frac{\partial}{\partial\tau}  \Big[\big(q & - & 1-\Delta_X
f\big)V  -  \tau T_h(V)\Big]  +  \Delta_X(Vf_{\tau})  = -(\Delta_X
f_{\tau})V+\big(q-1-\Delta_X f\big) \frac{\partial V}{\partial\tau}\\
& - &
T_h(V)-\tau\frac{\partial}{\partial\tau}\big(T_h(V)\big)+f_{\tau}\Delta_X
V+\nabla_0 V\cdot\nabla_0 f_{\tau}+V\Delta_Xf_{\tau}\\ & =
&\big(q-1-\Delta_X f\big) \frac{\partial
V}{\partial\tau}-\frac{\partial
V}{\partial\tau}-\nabla_0h\cdot\nabla_0V-\tau T_h\big(\frac{\partial
V}{\partial\tau}\big)-\tau\nabla_0 h_{\tau}\cdot\nabla_0V
\\ & +& f_{\tau}\Delta_X V+\nabla_0 V\cdot\nabla_0 f_{\tau}\\ & =
&\big(q-2-\Delta_X f\big) \frac{\partial
V}{\partial\tau}+\nabla_0V\cdot\Big(\nabla_0 f_{\tau}
-\nabla_0h-\tau\nabla_0h_{\tau}\Big)-\tau T_h\big(\frac{\partial
V}{\partial\tau}\big)+f_{\tau}\Delta_X V\\ & = & \big((q-2)-(\tau
T_h+\Delta_X f)\big) \frac{\partial
V}{\partial\tau}+f_{\tau}\Delta_X V.
\end{eqnarray*} Substituting the last expression into~\eqref{heq3}
we finally deduce
\begin{eqnarray}\label{heq4}(\Delta_X-\frac{\partial}{\partial u})\Big(\frac{e^{-\frac{f}{u}}(Vf_{\tau})}{u^{q}}
\Big) & = &
\Big(\frac{e^{-\frac{f}{u}}}{u^{q}}\Big)\Big(\big((q-2)-(\tau
T_h+\Delta_X f)\big) \frac{\partial
V}{\partial\tau}+f_{\tau}\Delta_X V\Big)\nonumber
\\ & + &\frac{1}{u^{q}}\frac{\partial}{\partial\tau}\Big(\big(e^{-\frac{f}{u}}\big)
\big(\tau T_h(V)+(\Delta_X f-q+1)V\big)\Big)\\ & + &
\frac{1}{u^{q+1}}\frac{\partial}{\partial\tau}\Big(e^{-\frac{f}{u}}\tau
f_{\tau}V\Big).\nonumber\end{eqnarray}

If we were suppose that $e^{-\frac{f}{u}}\tau f_{\tau}V$,
$e^{-\frac{f}{u}}\Big(\tau T_h(V)+(\Delta_X f-q+1)V\Big)$ tends to
$0$ as $\tau\to\pm\infty$, and the function $V$  satisfies the
transport equation
\begin{equation}\label{eq:transport}
\big((q-2)-(\tau T_h+\Delta_X f)\big) \frac{\partial
V}{\partial\tau}+f_{\tau}\Delta_X V=0,
\end{equation} then the heat kernel $P_u(\zeta_1,\zeta_2,\eta)$ would be found in the form~\eqref{eq:heatkernel}. Observe that we did not use any specific form of the auxiliary function $h$. It is sufficient  that $h$ satisfies the Hamilton--Jacobi equation. Take the modified action $g(\zeta_1,\zeta_2,\eta,\psi_1,-\psi_1,\tau)$ as the function $h(\tau)$ .

Let us present the coefficient $\Delta_X f$. Since $\frac{\partial
f}{\partial\zeta_1}=\psi_1$, $\frac{\partial
f}{\partial\zeta_2}=\psi_2$, $\frac{\partial
f}{\partial\eta}=\theta$, we have $\Delta_X f=2\frac{\partial^2
f}{\partial\eta^2}=2\frac{\partial\theta}{\partial\eta}$. For the
function $\theta$  the expression
\begin{equation}\label{eq:theta}\theta^2(\eta)=\theta_0^2-\tau^2\cotan^22\eta,\end{equation}
is valid, where $\theta_0$ can be found from the equation
\begin{equation}\label{eq:3}\cos
2\eta=\sqrt{\frac{\theta_0^2}{\tau^2+\theta_0^2}}\sin
8\sqrt{\tau^2+\theta_0^2},\end{equation} where we take into
consideration the sigh $(+)$ in the solution~\eqref{eq:eta}.
Differentiating~\eqref{eq:theta} one finds
$$\frac{\partial\theta(\eta)}{\partial\eta}=\frac{1}{2\theta(\eta)}\Big(
\frac{\partial\theta_0^2}{\partial\eta}+4\tau^2\frac{\cos
2\eta}{\sin^3 2\eta}\Big).$$ Expression~\eqref{eq:3} gives
$$\frac{\partial\theta_0^2}{\partial\eta}=
\frac{4(\tau^2+\theta_0^2)\sin 2\eta}{
\frac{\tau^2}{|\theta_0|\sqrt{\tau^2+\theta_0^2}}\sin
8\sqrt{\tau^2+\theta_0^2}+8|\theta_0|\cos
8\sqrt{\tau^2+\theta_0^2}}.$$ Observe that $$\frac{\sin
8\sqrt{\tau^2+\theta_0^2}}{\sqrt{\tau^2+\theta_0^2}}=-\frac{\cos
2\eta}{|\theta_0|},$$ and $$|\theta_0|\cos
8\sqrt{\tau^2+\theta_0^2}=\pm\cos2\eta\sqrt{\theta_0^2\tan^22\eta-\tau^2}.$$
This implies
$$\frac{\partial\theta_0^2}{\partial\eta}=\frac{4(\tau^2+\theta_0^2)\tan
2\eta}{\pm8\sqrt{\theta_0^2\tan^2
2\eta-\tau^2}-\frac{\tau^2}{\theta_0^2}},$$ which leads to the
sub-Laplacian  $$\Delta_X
f=\frac{1}{\theta(\eta)}\Big(4\tau^2\frac{\cos 2\eta}{\sin^3
2\eta}+\frac{4(\tau^2+\theta_0^2)\tan 2\eta}{\pm
8\sqrt{\theta_0^2\tan^2
2\eta-\tau^2}-\frac{\tau^2}{\theta_0^2}}\Big),$$ where
$\theta(\eta)$ is given by~\eqref{eq:theta} and $\theta_0$ can be
found from~\eqref{eq:3}. Since the coefficients of the transport
equation given by $\Delta_X f$ depend on $\eta$, we can not assume that 
$V$ depends only on $\tau$ contrasting with  the
case of   nilpotent groups. In fact, it depends also on $\eta$,
and we have to take into account the term
$\nabla_0g\cdot\nabla_0V$ in the transport equation. Since the vector fields $X$ and $Y$, and the functions $f(\zeta_1,\zeta_2,\eta,\tau,-\tau)$, and  $g(\zeta_1,\zeta_2,\eta,\psi_1,-\psi_1,\eta)$ depend only on $\eta$ and the difference $\zeta_1-\zeta_2$,  the coefficients of the transport equation depend only on $\eta$ and  $\zeta_1-\zeta_2$. We expect that the solution of the transport equation $V$ will be a function of $\tau$, $\eta$, and  $\zeta_1-\zeta_2$. This assumption makes the problem of finding volume element on $\mathbb S^3$ much more difficult, than in the case of nilpotent groups, where the solution of the transport equation was a function only on $\eta$.

\subsection{Volume element for the Green function} Following the geometric method elaborated by Beals, Gaveau, and Greiner, we suppose that the fundamental solution of the sub-Laplacian equation with the sub-Laplacian $\Delta_X=X^2+Y^2$, has the form \begin{equation}\label{green0} K(\zeta_i^0,\eta_0,\zeta_i,\eta)=\int_{+\infty}^{-\infty}\frac{w(t,\zeta_i,\eta)}{g(t,\zeta_i,\eta)}\,dt,\end{equation} where $w$ is the volume element over the characteristic variety at $(\zeta_i^0,\eta_0)$ which is just the straight line in our case. The function $g$ is the modified action defined in Section~\ref{Modif.function},  satisfying the Hamilton--Jacobi equation. In~\cite{BealsGaveauGreiner11} the authors elaborated the method of finding the volume element for the two step groups. Applying to our case, two step means that the matrix 
\begin{equation*}
-\frac{1}{2}\omega([X_i,X_j])=
\left[\array{rr}0 & -1 
\\
1 & 0 
\endarray\right],
\end{equation*} is non-singular at every point. Here $[X_i,X_j]$ are all possible commutators of $X$, $Y$, and $\omega$ is the one-form~\eqref{forma}. 
Here we give only a general scheme referring  the reader to~\cite{BealsGaveauGreiner11} for the details. 

Assume that the volume element $w$ is written as $w=W\mathcal E$, with $\E=-\frac{\partial g}{\partial t}=H(\nabla g)$, where we keep the notations of the previous subsections. In order $K$ to be the Green function, we  verify \begin{eqnarray}\label{green1}
0 & = & \int_{-\infty}^{+\infty}\Delta_X\Big(\frac{W\E}{g}\Big)\,dt\nonumber
\\ & = & \Big(\frac{\Delta_X(W\E)}{g}-\frac{\nabla_0g\cdot\nabla_0(W\E)}{g^2}+\frac{2W\E H(\nabla g)}{g^3}\Big)\,dt.
\end{eqnarray} Applying the Hamilton--Jacobi equation and integrating by parts twice, we continue~\eqref{green1} as follows

\begin{eqnarray*}
& = & \int_{-\infty}^{+\infty}\Big[\E \Delta_X(W)+\big(T_g+\Delta_X g\big)\frac{\partial W}{\partial t}\Big]\frac{dt}{g}\nonumber
\\
& + & \int_{-\infty}^{+\infty}\frac{\partial}{\partial t}\Big(\frac{W\E}{g^2}\Big)\,dt
\\
& - & \int_{-\infty}^{+\infty}\frac{\partial}{\partial t}\Big(\frac{T_g(W)+(\nabla_0 g)W}{g}\Big)\,dt.\nonumber
\end{eqnarray*}
Then, similarly to the case of the heat kernel, we assume that $g$, $\E=-\frac{\partial g}{\partial t}$, and $W$ are  suitable to satisfy $\frac{W\E}{g^2}\to 0$ and $\frac{T_h(W)+(\nabla_0 g)W}{g}\to 0$ as $t\to \pm\infty$. One comes to the conclusion that if $W$ satisfies the transport equation $$\E \Delta_X(W)+\big(T_g+\Delta_X g\big)\frac{\partial W}{\partial t}=0,$$ then the Green function is given by~\eqref{green0}.

\end{document}